\documentclass[a4paper,11pt]{amsart}
\usepackage{amsfonts,amssymb}
\author{Ivan V. Arzhantsev}
\thanks{This paper is to appear in the London Mathematical Society
Lecture Notes Series, published by Cambridge University Press,
who will own the copyright of the paper.
Supported by the London Mathematical Society, by RFBR grant
03-01-06252, CRDF grant RM1-2543-MO-03, and the RF
President grant MK-1279.2004.1}
\title{Affine embeddings of homogeneous spaces}
\date{March 22, 2005}

\newcommand{\kk}{\mathbb{K}}
\newcommand{\CC}{\mathbb{C}}

\newcommand{\PP}{\mathbb{P}}
\newcommand{\QQ}{\mathbb{Q}}
\newcommand{\HH}{\mathbb{H}}
\newcommand{\ZZ}{\mathbb{Z}}
\newcommand{\NN}{\mathbb{N}}
\newcommand{\kd}{{\rm Kdim\,}}
\newcommand{\modd}{{\mbox{mod}}}
\newcommand{\codim}{\mathop{\mathrm{codim\,}}}
\newcommand{\trdeg}{\mathop{\mathrm{tr.deg\,}}}
\newcommand{\tr}{\mathop{\mathrm{tr\,}}}
\newcommand{\cchar}{\mathop{\mathrm{char\,}}}
\newcommand{\Aut}{\mathop{\mathrm{Aut\,}}}
\newcommand{\End}{\mathop{\mathrm{End\,}}}
\newcommand{\Ker}{\mathop{\mathrm{Ker\,}}}
\newcommand{\Mat}{\mathop{\mathrm{Mat\,}}}
\newcommand{\Hom}{\mathop{\mathrm{Hom}}}
\newcommand{\Spin}{\mathop{\mathrm{Spin}}}
\newcommand{\Sp}{\mathop{\mathrm{Spec\,}}}
\newcommand{\rk}{\mathop{\mathrm{rk}}}
\newcommand{\A}{\mathfrak A}
\newcommand{\B}{\mathfrak B}
\newcommand{\C}{\mathfrak C}
\newcommand{\D}{\mathfrak D}
\newcommand{\I}{\mathfrak I}

\newtheorem{theorem}{Theorem}
\theoremstyle{conjecture}

\theoremstyle{definition}
\newtheorem{definition}{Definition}
\newtheorem{example}{Example}
\theoremstyle{proposition}
\newtheorem{proposition}{Proposition}
\theoremstyle{lemma}
\newtheorem{lemma}{Lemma}
\theoremstyle{corollary}
\newtheorem{corollary}{Corollary}
\theoremstyle{remark}
\newtheorem{remark}{Remark}
\begin{document}

\begin{abstract}

Let $G$ be a reductive algebraic group and $H$ a closed subgroup of $G$.
An affine embedding of the homogeneous space $G/H$
is an affine $G$-variety with an open $G$-orbit isomorphic to $G/H$.
The homogeneous space $G/H$ admits an affine embedding if and only if $G/H$
is a quasi-affine algebraic variety.

We start with some basic properties of
affine embeddings and consider the cases, where the theory is well-developed:
toric varieties, normal $SL(2)$-embeddings, $S$-varieties,
and algebraic monoids.
We discuss connections between the theory of affine embeddings and
Hilbert's 14th problem via a theorem of Grosshans. We characterize
affine homogeneous spaces $G/H$ such that any affine embedding of $G/H$
contains a finite number of $G$-orbits. The maximal value of modality over
all affine embeddings of a given affine homogeneous space $G/H$ is computed
and the group of equivariant automorphisms of an embedding is studied.

As applications of the theory, we describe
invariant algebras on homogeneous spaces
of a compact Lie group and $G$-algebras with finitely generated invariant
subalgebras.

\medskip

{\it AMS 2000 Math. Subject Classification: Primary 13A50, 14M17,
14R20; \ Secondary 14L30, 14M25, 14R05, 22C05, 32M12, 46J10}

\end{abstract}

\maketitle


\small

\tableofcontents

\normalsize

\section*{Introduction}

Throughout the paper
$G$ denotes a connected reductive algebraic
group, unless otherwise specified, and $H$ is an algebraic subgroup of $G$.
All groups and algebraic varieties are considered over an algebraically closed
field $\kk$ of characteristic zero, unless otherwise specified.
Let $\kk[X]$ be the algebra of regular functions on an algebraic variety $X$
and $\kk(X)$ the field of rational functions on $X$ provided $X$ is
irreducible. Our general references are:
\cite{hum} for algebraic groups, and \cite{pv2}, \cite{kr}, \cite{gr} for
algebraic transformation groups and invariant theory.

\smallskip

{\bf Affine embeddings: definitions.}\
Let us recall that an irreducible
algebraic $G$-variety $X$ is said to be an
\emph{embedding} of the homogeneous space $G/H$ if $X$ contains an open $G$-orbit
isomorphic to $G/H$. We
shall denote this by $G/H\hookrightarrow X$.
Let us say that
an embedding $G/H\hookrightarrow X$ is \emph{affine} if the variety
$X$ is affine.
In many problems of invariant theory, representation
theory, and other branches of mathematics, only affine
embeddings of homogeneous spaces arise. This is why it is reasonable
to study specific properties of affine embeddings in the framework of
a well-developed general embedding theory.

\smallskip

{\bf Which homogeneous spaces admit an affine embedding ?}\
It is easy to show that a homogeneous space $G/H$ admits
an affine embedding if and only if $G/H$ is quasi-affine (as an
algebraic variety).
In this situation, the subgroup
$H$ is said to be \emph{observable} in $G$.
A closed subgroup $H$ of $G$ is observable if and only if there
exist a rational finite-dimensional $G$-module $V$ and
a vector $v\in V$ such that the stabilizer $G_v$ coincides with $H$.
(This follows from the fact that any affine $G$-variety may be realized
as a closed invariant subvariety in a finite-dimensional
$G$-module~\cite[Th.1.5]{pv2}.)
There is a nice group-theoretic description of
observable subgroups due to A.~Sukhanov:
a subgroup $H$ is observable in $G$ if and only if
there exists a quasi-parabolic
subgroup $Q\subset G$ such that $H\subset Q$ and
the unipotent radical $H^u$
is contained in the unipotent radical $Q^u$,
see~\cite{su}, \cite[Th.7.3]{gr}.
(Let us recall that a subgroup $Q$ is said to be {\it quasi-parabolic}
if $Q$ is the stabilizer of a highest weight vector in some $G$-module $V$.)

It follows from Chevalley's theorem that any subgroup $H$
without non-trivial characters (in particular, any unipotent subgroup) is
observable.
By Matsushima's criterion, a homogeneous space
$G/H$ is affine if and only if $H$ is reductive.
(For a simple proof, see~\cite{lu2}; a characteristic-free proof can be found
in~\cite{ri}.) In particular, any reductive subgroup is observable.
A description of affine homogeneous spaces $G/H$ for non-reductive $G$
is still an open problem.

\smallskip

{\bf Complexity of reductive group actions.}\
Now we are going to define the notion of complexity,
which we shall encounter many times in the text.
Let us fix the notation. By $B=TU$ denote a Borel subgroup of $G$
with a maximal torus $T$ and the unipotent radical $U$. By definition, the
\emph{complexity} $c(X)$ of a $G$-variety $X$ is
the codimension of a $B$-orbit of general position in $X$
for the restricted action $B:X$. This notion firstly appeared in
\cite{lv} and~\cite{vi1}. Now it plays a central role in embedding theory.
By Rosenlicht's theorem, $c(X)$ is equal to the transcendence
degree of the field $\kk(X)^B$ of rational $B$-invariant
functions on~$X$. A normal $G$-variety $X$ is called
\emph{spherical} if $c(X)=0$, or, equivalently, $\kk(X)^B=\kk$.
A homogeneous space $G/H$ and a subgroup $H\subseteq G$ are said
to be {\it spherical} if $G/H$ is a spherical $G$-variety.

\smallskip

{\bf Rational representations, the isotypic decomposition, and $G$-algebras.}\
A linear action of $G$ in vector space $W$ is said to be {\it rational}
if for any vector $w\in W$ the linear span
$\langle Gw \rangle$ is finite-dimensional
and the action $G:\langle Gw\rangle$ defines a representation of an algebraic
group. Since any finite-dimensional representation
of $G$ is completely reducible, it is easy to prove that $W$ is a direct sum
of finite-dimensional simple $G$-modules. Let $\Xi_+(G)$ be the semigroup
of dominant weights of $G$. For any $\lambda\in\Xi_+(G)$, by $W_{\lambda}$
denote
the sum of all simple submodules in $W$ of highest weight $\lambda$.
The subspace $W_{\lambda}$ is called
{\it an isotypic component} of $W$ of weight
$\lambda$, and the decomposition
$$
  W=\oplus_{\lambda\in\Xi_+(G)} W_{\lambda}
$$
is called {\it the isotypic decomposition} of $W$.

If $G$ acts on an affine variety $X$, the linear action $G:\kk[X]$,
$(gf)(x):=f(g^{-1}x)$, is rational~\cite[Lemma~1.4]{pv2}.
(Note that for irreducible $X$ the action on rational functions $G:\kk(X)$
defined by the same formula is not rational.) The isotypic decomposition
$$
  \kk[X]=\oplus_{\lambda\in\Xi_+(G)} \kk[X]_{\lambda}
$$
and its interaction with the multiplicative structure on $\kk[X]$ give
important technical tools for the study of affine embeddings.

An affine $G$-variety $X$ is spherical if and only if $\kk[X]_{\lambda}$
is a simple $G$-module for any $\lambda\in\Xi_+(G)$~\cite{vk}.

Suppose that $\A$ is a commutative associative algebra
with unit over $\kk$. If $G$ acts on $\A$ by automorphisms and
the action $G:\A$ is rational, we say that $\A$ is a {\it $G$-algebra}.
The algebra $\kk[X]$ is a $G$-algebra for any affine $G$-variety $X$.
Moreover, any finitely generated
$G$-algebra without nilpotents arises in this way.

\medskip

We conclude Introduction by a review of the contents of this survey.

One of the pioneer works in embedding theory was a  classification
of normal affine $SL(2)$-embeddings due to V.~L.~Popov,
see~\cite{pop1}, \cite{kr}.
At the same period (early seventies),
the theory of toric varieties was developed. A toric variety may be considered
as an equivariant embedding of an algebraic torus $T$.
Such embeddings are described in terms
of convex fans. Any cone in the fan of a toric variety $X$ represents
an affine toric variety.
This reflects the fact that $X$ has a covering
by $T$-invariant affine charts.
In 1972, V.~L.~Popov and E.~B.~Vinberg~\cite{vp1} described affine embeddings
of quasi-affine homogeneous spaces $G/H$, where $H$ contains a maximal unipotent
subgroup of $G$.  In Section~\ref{Sec1}, we discuss briefly these results
together with the more recent remarkable classification of algebraic
monoids with a reductive group $G$ as the group of invertible elements
(E.~B.~Vinberg~\cite{vi2}), which is nothing else
but the classification of affine
embeddings of the space $(G\times G)/\Delta(G)$, where $\Delta(G)$
is the diagonal subgroup.

In Section~\ref{Sec2} we consider connections of the theory of affine
embeddings with Hilbert's 14th problem. Let $H$ be an observable subgroup
of $G$. By the Grosshans theorem, the
following conditions are equivalent: 1) the algebra of invariants $\kk[V]^H$
is finitely generated for any $G$-module $V$;
\ 2) the algebra of regular functions $\kk[G/H]$ is finitely generated;
\ 3) there exists a (normal) affine embedding $G/H\hookrightarrow X$
such that $\codim_X (X\setminus (G/H))\ge 2$ (such an embedding is called
{\it the canonical embedding} of $G/H$).

It was proved by F.~Knop that if $c(G/H)\le 1$ then the algebra $\kk[G/H]$
is finitely generated. This result provides a large class of subgroups
with a positive solution of Hilbert's 14th problem. In particular,
Knop's theorem together with Grosshans' theorem
on the unipotent radical $P^u$ of a parabolic subgroup $P\subset G$
cover almost all known results on Popov-Pommerening's
conjecture (see~\ref{subsec2.2}). We study the canonical embedding of $G/P^u$
from a geometric view-point.
Finally, we mention counterexamples to Hilbert's 14th problem due to
M.~Nagata, P.~Roberts, and R.~Steinberg.

\smallskip

In Section~\ref{Sec3}, we introduce the notion of {\it an affinely closed
space}, i.e., an affine homogeneous space admitting no non-trivial affine
embeddings, and discuss the result of D.~Luna related to this notion.
(We say that an affine embedding $G/H\hookrightarrow X$ is {\it trivial}
if $X=G/H$.) Affinely closed spaces of an arbitrary affine algebraic group
are characterized and
some elementary properties of affine embeddings are formulated.

\smallskip

Section~\ref{Sec4} is devoted to affine embeddings with a finite number of
orbits. We give a characterization of affine homogeneous spaces $G/H$ such
that any affine embedding of $G/H$ contains a finite number of orbits.
More generally, we compute the maximal number of parameters
in a continuous family of $G$-orbits over all affine embeddings of a given
affine homogeneous space $G/H$. The group of equivariant automorphisms
of an affine embedding is also studied here.

\smallskip

Some applications of the theory of affine embeddings to functional analysis are
given in Section~\ref{Sec5}.
Let $M=K/L$ be a homogeneous space of a connected compact Lie group $K$, and
$C(M)$ the commutative Banach algebra of all complex-valued continuous functions
on $M$. The $K$-action on $C(M)$ is defined by the formula $(kf)(x)=f(k^{-1}x), \
k\in K, \ x\in M$. We shall say that $A$ is an {\it invariant algebra} on $M$
if $A$ a $K$-invariant uniformly closed subalgebra with unit in $C(M)$.
Denote by $G$ (resp. $H$) the complexification
of $K$ (resp. $L$). Then $G$ is a reductive algebraic group with a reductive
subgroup $H$. There exists a correspondence between
finitely generated invariant algebras
on $M$ and affine embeddings of $G/F$ with some additional data,
where $F$ is an observable subgroup of $G$ containing $H$. This correspondence
was introduced by V.~M.~Gichev~\cite{gich},
I.~A.~Latypov~\cite{lat},~\cite{latdis}
and, in a more algebraic way, by E.~B.~Vinberg.
We give a precise formulation of this correspondence
and reformulate some results on affine embeddings in terms of invariant
algebras. Some results of this section are new and not published elsewhere.

\smallskip

The last section is devoted to $G$-algebras.
It is easy to prove that any subalgebra in the polynomial
algebra $\kk[x]$ is finitely generated. On the other hand, one can construct
many non-finitely generated subalgebras in $\kk[x_1,\dots,x_n]$ for $n\ge 2$.
More generally, every subalgebra in an associative commutative finitely generated
integral domain $\A$ with unit is finitely generated if and only if
$\kd \A\le 1$, where $\kd \A$ is Krull dimension of $\A$ (Proposition~\ref{ga}).
In Section~\ref{secap2} we obtain an equivariant version of this result.
The problem was motivated by the study of invariant algebras in the previous
section.
The proof of the main result (Theorem~\ref{tmain})
is based on a geometric method for constructing a non-finitely
generated subalgebra in a finitely generated
$G$-algebra and on properties of affine embeddings
obtained above. In particular, the notion of an affinely closed space is crucial
for the classification of $G$-algebras with finitely
generated invariant subalgebras.
The arguments used in this text are slightly different
from the original ones~\cite{ar5}. A characterization of $G$-algebras with
finitely generated invariant subalgebras for non-reductive $G$ is also given
in this section.

\medskip

{\bf Acknowledgements.}\ These notes were initiated by my visit to
the Manchester University in March, 2003. I am grateful to this
institution for hospitality, to Prof. A.~Premet for invitation and
organization of this visit, and to the London Mathematical Society
for financial support. The work was continued during my stay at
Institut Fourier (Grenoble) in April-July, 2003. I would like to
express my gratitude to this institution and especially to Prof.
M.~Brion for the invitation, and for numerous remarks and
suggestions. Thanks are also due to F.~D.~Grosshans and
D.~A.~Timashev for useful remarks.


\section{Remarkable classes of affine embeddings}\label{Sec1}

\subsection{Affine toric varieties}

We begin with some notation. Let $T$ be an algebraic torus and $\Xi(T)$ the
lattice of its characters. A $T$-action on an affine variety $X$ defines
a $\Xi(T)$-grading on the algebra $\kk[X]=\oplus_{\chi\in \Xi(T)} \kk[X]_{\chi}$,
where $\kk[X]_{\chi}=\{ f \mid tf=\chi(t)f \ {\rm for\ any} \ t\in T\}$.
(This grading is nothing else but the isotypic decomposition, see Introduction.)
If $X$ is irreducible, then the set $L(X)=\{ \chi \mid \kk[X]_{\chi}\ne 0\}$
is a submonoid in $\Xi(T)$.

\begin{definition}
An affine toric variety $X$ is a normal affine $T$-variety with an open
$T$-orbit isomorphic to $T$.
\end{definition}

Below we list some basic properties of $T$-actions:

\smallskip

$\bullet$
An action $T:X$ has an open orbit if and only if $\dim\kk[X]_{\chi}=1$ for
any $\chi\in L(X)$. In this situation $\kk[X]$ is $T$-isomorphic to the
semigroup algebra $\kk L(X)$.

\smallskip

$\bullet$ An action $T:X$ is effective if and only if the subgroup in $\Xi(T)$
generated by $L(X)$ coincides with $\Xi(T)$.

\smallskip

$\bullet$ Suppose that $T:X$ is an effective action with an open orbit.
Then the following conditions are equivalent:

1) $X$ is normal;

2) the semigroup algebra $\kk L(X)$ is integrally closed;

3) if $\chi\in \Xi(T)$ and there exists $n\in\NN$, $n>0$,
such that $n\chi\in L(X)$,
then $\chi\in L(X)$ \ (the saturation condition);

4) there exists a solid convex polyhedral cone $K$
in $\Xi(T)\otimes_{\ZZ}\QQ$ such that $L(X)=K\cap \Xi(T)$.

In this situation, any $T$-invariant radical ideal of $\kk[X]$
corresponds to the subsemigroup $L(X)\setminus M$ for a fixed face $M$ of the
cone $K$. This correspondence defines a bijection between $T$-invariant
radical ideals of $\kk[X]$ and faces of $K$.

\medskip

The proof of these properties can be found for example, in~\cite{fu}.
Summarizing all the results, we obtain

\begin{theorem}\label{teotor}
1)\ Affine toric varieties are in one-to-one correspondence with
solid convex polyhedral cones in the space $\Xi(T)\otimes_{\ZZ}\QQ$;

2) $T$-orbits on a toric variety are in one-to-one correspondence with
faces of the cone.
\end{theorem}

The classification of affine toric varieties will serve us as a sample
for studing more complicated classes of affine embeddings.
Generalizations of a combinatorial description of toric
varieties were obtained for spherical varieties~\cite{lv}, \cite{kn1},
\cite{br2}, and for embeddings of complexity one~\cite{tim}. In this
more general context, the idea that normal $G$-varieties may be described
by some convex cones becomes rigorous through the method of $U$-invariants
developed by D.~Luna and Th.~Vust. The essence of this method is contained
in the following theorem (see~\cite{vu}, \cite{kr}, \cite{pc}, \cite{gr}).

\begin{theorem}\label{teouinv}
Let $\A$ be a $G$-algebra
and $U$ a maximal unipotent subgroup of $G$.
Consider the following properties of an algebra:

1) it is finitely generated;

2) it has no nilpotent elements;

3) it has no zero divisors;

4) it is integrally closed.

If (P) is any of these properties, then the algebra $\A$ has property (P)
if and only if the algebra $\A^U$ has property (P).
\end{theorem}

We try to demonstrate briefly some applications
of the method of $U$-invariants in the following subsections.


\subsection{Normal affine $SL(2)$-embeddings}

Suppose that the group $SL(2)$ acts on a normal affine
variety $X$ and there is a
point $x\in X$ such that the stabilizer of $x$ is trivial and the orbit
$SL(2)x$ is open in $X$. We say in this case that $X$ is a {\it normal
$SL(2)$-embedding}.

Let $U_m$ be a finite extension of the standard maximal unipotent subgroup in $SL(2)$:
$$
  U_m=\left\{ \left( \begin{array}{cc}
                      \epsilon & a \\
                      0 &\epsilon^{-1}
                      \end{array} \right)
\left| \right. \ \epsilon^m=1, \ a\in \kk \right\}.
$$

\begin{theorem}~\cite{pop1} \label{teopop}
Normal non-trivial $SL(2)$-embeddings are in one-to-one
correspondence with rational numbers $h\in (0,1]$. Furthermore,

1) $h=1$ corresponds to a (unique)
smooth $SL(2)$-embedding with two orbits: $X=SL(2)\cup SL(2)/T$;

2) if $h=\frac{p}{q}<1$ and $(p,q)=1$, then
$X=SL(2)\cup SL(2)/U_{p+q} \cup \{pt\}$, and $\{pt\}$ is an isolated singular
point in $X$.
\end{theorem}

The proof of Theorem~\ref{teopop} can be found in~\cite{pop1},
\cite[Ch.~3]{kr}. Here we give only some examples and explain what
the number $h$ (which is called the {\it height} of $X$) means in
terms of the algebra $\kk[X]$.

\begin{example}
1) The group $SL(2)$ acts tautologically on space $\kk^2$ and by
conjugation on space $\Mat(2\times 2)$. Consider the point
$$
x=\left\{ \left( \begin{array}{cc}
                  1 & 0 \\
                  0 & -1
                  \end{array} \right),
                  \left( \begin{array}{c}
                          1 \\
                          0
                         \end{array} \right) \right\}
                   \in \Mat(2\times 2)\times \kk^2
$$
and its orbit
$$
SL(2)x=\{(A,v) \mid \det A=-1, \tr A=0, Av=v, v\ne 0 \}.
$$
It is easy to see that the closure
$$
X=\overline{SL(2)x}=\{(A, v) \mid \det A=-1, \tr A=0, Av=v\}
$$
is a smooth $SL(2)$-embedding with two orbits, hence $X$ corresponds to ${h=1}$.

2) Let $V_d=\langle x^d, x^{d-1}y, \dots, y^d\rangle$
be the $SL(2)$-module of binary
forms of degree $d$. It is possible to check that
$$
  X=\overline{SL(2)(x, x^2y)}\subset V_1\oplus V_3
$$
is a normal $SL(2)$-embedding with the orbit decomposition
$X=SL(2)\cup SL(2)/U_3\cup \{pt\}$, hence $X$ corresponds to $h=\frac{1}{2}$.
\end{example}

An embedding $SL(2)\hookrightarrow X$, $g\to gx$
determines the injective homomorphism
$\A=\kk[X]\to \kk[SL(2)]$ with $Q\A=Q\kk[SL(2)]$,
where $Q\A$ is the quotient field of $\A$.
Let $U^-$ be the unipotent subgroup of $SL(2)$ opposite
to $U$. Then
$$
  \kk[SL(2)]^{U^-}=\{f\in\kk[SL(2)] \mid f(ug)=f(g), \ g\in SL(2), u\in U^- \}=
\kk[A,B],
$$
where $A\left( \begin{array}{cc}
               a & b \\
               c & d
              \end{array} \right)=a$ and
$B\left( \begin{array}{cc}
               a & b \\
               c & d
              \end{array} \right)=b$.

Below we list some facts (\cite[Ch.~3]{kr}) that allow to
introduce the height of an $SL(2)$-embedding $X$.

\smallskip

$\bullet$ If $\C$ is an integral $F$-domain, where $F$ is a
unipotent group, then $Q(\C^F)=(Q\C)^F$. In particular, if
$\C\subseteq \A$ and $Q\A=Q\C$, then
$Q(\A^{U^-})=Q(\C^{U^-})$.

\smallskip

$\bullet$ Suppose that $\lim_{t \to 0} \left( \begin{array}{cc}
                                 t & 0 \\
                                 0 & t^{-1}
                                 \end{array} \right)x$ exists.
Then $A\in\kk[SL(2)]\subset\kk(X)$ is regular on $X$.

\smallskip

$\bullet$ Let $\D\subset \kk[x,y]$ be a homogeneous integrally closed
subalgebra in the polynomial algebra such that $Q\D=\kk(x,y)$ and
$x\in \D$. Then $\D$ is generated by monomials.

In our situation, the algebra $\D=\A^{U^-}\subset \kk[A,B]$
is homogeneous because it is $T$-stable (since $T$ normalizes $U^-$).

\smallskip

$\bullet$ There exists rational $h\in (0,1]$ such that
$$
\A^{U^-}=\A(h)=\langle A^iB^j \mid \frac{j}{i}\le h \rangle.
$$
Moreover, for any rational $h\in (0,1]$ the subspace
$\langle SL(2)\A(h) \rangle\subset \kk[SL(2)]$ is a subalgebra.

\begin{remark}
While normal $SL(2)$-embeddings are parametrized by a discrete parameter
$h$, there are families of non-isomorphic non-normal $SL(2)$-embeddings
over a base of arbitrary dimension~\cite{ba}.
\end{remark}

\begin{remark}
A classification of $SL(2)$-actions on normal three-dimensional affine varieties
without open orbit can be found in~\cite{ar1},~\cite{ar2}.
\end{remark}


\subsection{$HV$-varieties and $S$-varieties}\label{sub1.3}

In this subsection we discuss the results of V.~L.~Popov and
E.~B.~Vinberg~\cite{vp1}. Suppose that $G$ is a connected and simply
connected semisimple group.

\begin{definition}
An {\it $HV$-variety} $X$ is the closure of the orbit of a highest
weight vector in a simple $G$-module.
\end{definition}

Let $V(\lambda)$ be the simple $G$-module with highest weight $\lambda$ and
$v_{\lambda}$ a highest weight vector in $V(\lambda)$.
Denote by $\lambda^*$ the
highest weight of the dual $G$-module $V(\lambda)^*$.

\smallskip

$\bullet$ \ $X(\lambda)=\overline{Gv_{\lambda^*}}$ is a normal
affine variety consisting
of two orbits: \ $X(\lambda)=Gv_{\lambda^*}\cup \{0\}$.

\smallskip

$\bullet$ \ $\kk[X(\lambda)]=\kk[Gv_{\lambda^*}]=\oplus_{m\ge 0}
\kk[X(\lambda)]_{m\lambda}$,
any isotypic component $\kk[X(\lambda)]_{m\lambda}$ is a simple $G$-module,
and
$\kk[X(\lambda)]_{m_1\lambda}\kk[X(\lambda)]_{m_2\lambda}=
\kk[X(\lambda)]_{(m_1+m_2)\lambda}$.

\smallskip

$\bullet$ The algebra $\kk[X(\lambda)]$ is a unique factorization domain
if and only if $\lambda$ is a
fundamental weight of $G$.

\begin{example}
1) The quadratic cone $KQ_n=\{ x\in\kk^n \mid x_1^2+\dots+x_n^2=0\}$
($n\ge 3$) is an
$HV$-variety for the tautological representation $SO(n):\kk^n$. (In fact, the
group $SO(n)$ is not simply connected and we consider the corresponding
module as a $\Spin(n)$-module.) It follows that $KQ_n$ is normal and
it is factorial if and only if $n\ge 5$.

2) The Grassmannian cone $KG_{n,m}$ ($n\ge 2,\ 1\le m\le n-1$) (i.e., the cone
over the projective variety of $m$-subspaces in $\kk^n$)
is an $HV$-variety
associated with the fundamental $SL(n)$-representation in space
$\bigwedge^m\kk^n$, hence it is factorial.
\end{example}

\begin{definition}
An irreducible affine variety $X$ with an action of a connected reductive
group $G$ is said to be {\it an $S$-variety} if $X$ has an open $G$-orbit
and the stabilizer of a point in this orbit contains a maximal unipotent
subgroup of $G$.
\end{definition}

Any $S$-variety may be realized as $X=\overline{Gv}$, where
$v=v_{\lambda_1^*}+\dots+v_{\lambda_k^*}$ is a sum of highest weight vectors
$v_{\lambda_i^*}$ in some $G$-module $V$. We have the isotypic decomposition
$$
 \kk[X]=\oplus_{\lambda\in L(X)} \kk[X]_{\lambda},
$$
where $L(X)$ is the semigroup generated by $\lambda_1,\dots,\lambda_k$,
any $\kk[X]_{\lambda}$ is a simple $G$-module, and
$\kk[X]_{\lambda}\kk[X]_{\mu}=\kk[X]_{\lambda+\mu}$.
The last condition determines uniquely (up to $G$-isomorphism)
the multiplicative structure on the $G$-module $\kk[X]$.
This shows that there is a bijection between $S$-varieties and finitely
generated submonoids in $\Xi_+(G)$.

Consider the cone $K=\QQ_+L(X)$. As in the toric case, normality of $X$ is
equivalent to the saturation condition for the semigroup $L(X)$, and
$G$-orbits on $X$ are in one-to-one correspondence with faces of $K$.
On the other hand,
there are phenomena which are specific for $S$-varieties. For example,
the complement to the open orbit in $X$ has codimension $\ge 2$ if and only if
$\ZZ L(X)\cap \Xi_+(G)\subseteq\QQ_+L(X)$
(this is never the case for non-trivial toric varieties).
For semisimple simply connected $G$,
an $S$-variety $X$ is factorial if and only if $L(X)$
is generated by fundamental weights.

Finally, we mention one more result on this subject. Say that an action
$G:X$ on an affine variety $X$
is {\it special} (or {\it horospherical})
if there is an open dense subset $W\subset X$ such that
the stabilizer of any point of $W$ contains a maximal unipotent subgroup of $G$.

\begin{theorem}~\cite{pc} \label{teospec}
The following conditions are equivalent:

1) the action $G:X$ is special;

2) the stabilizer of any point on $X$ contains a maximal unipotent subgroup;

3) $\kk[X]_{\lambda}\kk[X]_{\mu}\subseteq \kk[X]_{\lambda+\mu}$ for any
$\lambda, \mu\in \Xi_+(G)$.
\end{theorem}


\subsection{Algebraic monoids} The general theory of algebraic semigroups was
developed by M.~S.~Putcha, L.~Renner and E.~B.~Vinberg. In this subsection
we recall briefly the classification results following~\cite{vi2}.

\begin{definition}
An (affine) algebraic semigroup is an (affine) algebraic variety $S$ with
an associative multiplication
$$
 \mu:\ S\times S\to S,
$$
which is a morphism of algebraic varieties. An algebraic semigroup $S$ is
{\it normal} if $S$ is a normal algebraic variety.
\end{definition}

Any algebraic group is an algebraic semigroup. Another example
is the semigroup $\End(V)$ of endomorphisms of finite-dimensional vector
space $V$.

\begin{lemma}
An affine algebraic semigroup is isomorphic to a closed subsemigroup of $\End(V)$
for a suitable $V$. If $S$ has a unit, one may assume that it corresponds
to the identity map of $V$.
\end{lemma}

{\bf Proof.}\ The morphism $\mu: S\times S\to S$ induces the homomorphism
$\mu^*:\kk[S]\to\kk[S]\otimes\kk[S]$, $f(s) \mapsto F(s_1,s_2):=f(s_1s_2)$.
Hence $f(s_1s_2)=\sum_{i=1}^n f_i(s_1)h_i(s_2)$. Consider the linear action
$S:\kk[S]$ defined by $(s*f)(x)=f(xs)$.
One has $\langle Sf \rangle \subseteq \langle f_1,\dots,f_n \rangle$,
i.e., the linear span of any "$S$-orbit" in $\kk[S]$ is finite-dimensional and
the linear action $S: \langle Sf \rangle$ defines an algebraic representation
of $S$. Take as $V$ any finite-dimensional
$S$-invariant subspace of $\kk[S]$ containing a system of
generators of $\kk[S]$.

Suppose that $S$ is a monoid, i.e., a semigroup with unit.
We claim that the action $S:V$ defines a
closed embedding $\phi:S\to\End(V)$. Indeed, there are
$\alpha_{ij}\in\kk[S]$ such that $s*f_i=\sum_j \alpha_{ij}(s) f_j$.
The equalities $f_i(s)=(s*f_i)(e)=\sum_j\alpha_{ij}(s)f_j(e)$ show that
the homomorphism $\phi^*:\kk[\End(V)]\to\kk[S]$ is surjective.

The general case can be reduced to the previous one as follows: to any
semigroup $S$ one may add an element $e$ with relations $e^2=e$ and
$es=se=s$ for any $s\in S$. Then $\tilde S=S\sqcup \{e\}$ is an algebraic
monoid. $\diamondsuit$

\smallskip

If $S\subseteq\End(V)$ is a monoid, then any invertible element
of $S$ corresponds to an element of $GL(V)$.
Conversely, if the image of $s$ is invertible in $\End(V)$, then it is
invertible in $S$. Indeed, the sequence of closed subsets
$S\supseteq sS\supseteq s^2S\supseteq s^3S\supseteq\dots$ stabilizes,
and $s^kS=s^{k+1}S$ implies $S=sS$.
Hence the group $G(S)$ of invertible elements is open in $S$ and is an
algebraic group. Suppose that $G(S)$ is dense in $S$. Then
$S$ may be considered as an affine embedding
of $G(S)/\{e\}$ (with respect to left multiplication).

\begin{proposition}\label{prmon}
Let $G$ be an algebraic group. An affine embedding
$G/\{e\}\hookrightarrow S$ has a structure of an algebraic monoid with $G$
as the group of invertible elements if and only if the $G$-equivariant
$G$-action
on the open orbit by right multiplication can be extended to $S$, or,
equivalently, $S$ is an affine embedding of $(G\times G)/\Delta(G)$,
where $\Delta(G)$ is the diagonal in $G\times G$.
\end{proposition}

{\bf Proof.}\ If $S$ is an algebraic monoid with $G(S)=G$ and $G(S)$ is dense
in $S$, then $G\times G$
acts on $S$ by $((g_1,g_2), s) \mapsto g_1sg_2^{-1}$
and the dense open $G\times G$-orbit
in $S$ is isomorphic to $(G\times G)/\Delta(G)$.

For the converse, we give two independent proofs following
the chronological succession.

\smallskip

{\it Proof One (the reductive case).} (E.B.Vinberg~\cite{vi2})
An algebraic monoid $S$ is {\it reductive} if the group $G(S)$ is reductive
and dense in $S$.
The multiplication $\mu:G\times G\to G$ corresponds to
the comultiplication $\mu^*:\kk[G]\to \kk[G]\otimes\kk[G]$. Any
$(G\times G)$-isotypic component in $\kk[G]$ is a simple $(G\times G)$-module
isomorphic to $V(\lambda)^*\otimes V(\lambda)$
for $\lambda\in\Xi_+(G)$~\cite{kr}.
It coincides with the linear span of the matrix entries
of the $G$-module $V(\lambda)$. This shows that $\mu^*$ maps an
isotypic component to its tensor square, and for
any $(G\times G)$-invariant subspace $W\subset \kk[G]$ one has
$\mu^*(W)\subset W\otimes W$. Thus the spectrum $S$ of any
$(G\times G)$-invariant finitely generated
subalgebra in $\kk[G]$ carries the structure of an algebraic semigroup.
If the open $(G\times G)$-orbit in $S$ is isomorphic to
$(G\times G)/\Delta(G)$,
then $G(S)=G$. Indeed, $G$ is dense in $S$ and for any $s\in G(S)$
the intersection $sG\cap G\ne \emptyset$, hence $s\in G$.

\smallskip

{\it Proof Two (the general case).} (A.Rittatore~\cite{rit})
If the multiplication $\mu:G\times G\to G$ extends to a morphism
$\mu:S\times S\to S$, then $\mu$ is a multiplication because $\mu$ is associative
on $G\times G$. It is clear that $1\in G$ satisfies $1s=s1=s$ for all $s\in S$.
Consider the right and left actions of $G$ given by
$$
  G\times S\to S , \ \ \ gs=(g,1)s,
$$
$$
  S\times G\to S,  \ \ \ sg=(1,g^{-1}s).
$$ These actions define coactions $\kk[S]\to\kk[G]\otimes\kk[S]$
and $\kk[S]\to\kk[S]\otimes\kk[G]$, which are the restrictions to
$\kk[S]$ of the comultiplication $\kk[G]\to\kk[G]\otimes\kk[G]$.
Hence the image of $\kk[S]$ lies in
$$
  (\kk[G]\otimes\kk[S])\cap(\kk[S]\otimes\kk[G])=\kk[S]\otimes\kk[S],
$$
and we have a multiplication on $S$. The equality $G(S)=G$ may be proved
as above. $\diamondsuit$

\smallskip

Below we assume that $G$ is reductive.
For $\lambda_1, \lambda_2\in\Xi_+(G)$, by $\Xi(\lambda_1,\lambda_2)$ denote
the set of $\lambda\in\Xi_+(G)$ such that the $G$-module
$V(\lambda_1)\otimes V(\lambda_2)$
contains a submodule isomorphic to $V(\lambda)$. Since any
$(G\times G)$-isotypic component $\kk[G]_{(\lambda^*,\lambda)}$ in $\kk[G]$ is
the linear span of the matrix entries corresponding to the representation
of $G$ in $V(\lambda)$, the product
$\kk[G]_{(\lambda^*_1,\lambda_1)}\kk[G]_{(\lambda_2^*,\lambda_2)}$
is the linear span of the matrix entries corresponding to
$V(\lambda_1)\otimes V(\lambda_2)$. This shows that
$$
  \kk[G]_{(\lambda^*_1,\lambda_1)}\kk[G]_{(\lambda_2^*,\lambda_2)}=
\oplus_{\lambda\in\Xi(\lambda_1,\lambda_2)} \kk[G]_{(\lambda^*,\lambda)}.
$$

Since every $(G\times G)$-isotypic component in $\kk[G]$ is simple,
any $(G\times G)$-invariant
subalgebra in $\kk[G]$ is determined by the semigroup of dominant weights
that appear in its isotypic decomposition, and it is natural to classify
reductive algebraic monoids $S$ with $G(S)=G$ in terms of the semigroup
that determines $\kk[S]$ in $\kk[G]$.

\begin{definition}
A subsemigroup $L\subset\Xi_+(G)$ is {\it perfect} if it contains zero
and $\lambda_1,\lambda_2\in L$ implies $\Xi(\lambda_1,\lambda_2)\subset L$.
\end{definition}

\smallskip

Let $\ZZ\Xi_+(G)$ be the group generated by the semigroup $\Xi_+(G)$.
This group may be realized as the group of characters $\Xi(T)$
of a maximal torus of $G$.

\smallskip

\begin{theorem}~\cite{vi2}
A subset $L\subset\Xi_+(G)$ defines an affine algebraic monoid $S$ with
$G(S)=G$ if and only if $L$ is a perfect finitely generated subsemigroup
generating the group $\ZZ\Xi_+(G)$.
\end{theorem}

The classification of normal affine reductive monoids is more constructive.
Fix some notation. The group $G=ZG'$ is an almost direct product of
its center $Z$ and the derived subgroup $G'$.
Fix a Borel subgroup $B_0$ and a maximal torus $T_0\subset B_0$ in $G'$.
Then $B=ZB_0$ (resp. $T=ZT_0$) is a Borel subgroup (resp. a maximal torus)
in $G$. By $N$ (resp. $N_0$, $N_1$) denote $\QQ$-vector space
$\Xi(T)\otimes_{\ZZ}\QQ$ (resp. $\Xi(T_0)\otimes_{\ZZ}\QQ$,
$\Xi(Z)\otimes_{\ZZ}\QQ$). Then $N=N_1\oplus N_0$.
The semigroup of dominant weights $\Xi_+(G)$ (with respect to $B$) is
a subsemigroup in $\Xi(T)\subset N$. By
$\alpha_1,\dots,\alpha_k\in N_1$ denote the simple roots of $G$ with respect
to $B$, and by $C\subset N$ (resp. $C_0\subset N_0$) the positive Weyl
chamber for the group $G$ (resp. $G'$)
with respect to $\alpha_1,\dots,\alpha_k$.

\begin{theorem} \cite{vi2}
A subset $L\subset\Xi_+(G)$ defines a normal affine
algebraic monoid $S$ with $G(S)=G$
if and only if $L=\Xi_+(G)\cap K$, where $K$ is a closed convex polyhedral
cone in $N$ satisfying the conditions:

1) $-\alpha_1,\dots,-\alpha_k\in K$;

2) the cone $K\cap C$ generates $N$.

\noindent The monoid $S$ has a zero if and only if:

3) the cone $K\cap N_1$ is pointed;

4) $K\cap C_0=\{ 0\}$.
\end{theorem}

A characteristic-free approach to the classification of reductive
algebraic monoids via the theory of spherical varieties was developed
in~\cite{rit}. Another interesting result of~\cite{rit} is that any
reductive algebraic monoid is affine. Recently A.~Rittatore announced
a proof of the fact that any algebraic monoid with an affine algebraic
group of invertible elements is affine.


\section{Connections with Hilbert's 14th Problem}\label{Sec2}

\subsection{Grosshans subgroups and the canonical embedding} \
Let $H$ be a closed subgroup of $GL(V)$. Hilbert's 14th problem
(in its modern version) may be formulated as follows: characterize
subgroups $H$ such that the algebra of polynomial
invariants $\kk[V]^H$ is finitely
generated. It is a classical result that for reductive $H$ the algebra $\kk[V]^H$
is finitely generated. For non-reductive linear groups this problem
seems to be very far from a complete solution.

\smallskip

\begin{remark}
Hilbert's original statement of the problem was the following:

\smallskip

{\it For a field $\kk$, let $\kk[x_1,\dots,x_n]$ denote the polynomial ring
in $n$ variables over $\kk$, and let $\kk(x_1,\dots,x_n)$ denote its
field of fractions. If $K$ is a subfield of $\kk(x_1,\dots,x_n)$ containing
$\kk$, is $K\cap\kk[x_1,\dots,x_n]$ finitely generated over $\kk$ ? }

\smallskip

Since $\kk[V]^H=\kk[V]\cap\kk(V)^H$, our situation may be regarded as
a particular case of the general one.
\end{remark}

\smallskip

Let us assume that $H$ is a subgroup of a bigger reductive group $G$
acting on $V$. (For example, one may take $G=GL(V)$.) The intersection of a
family of observable subgroups in $G$ is an observable subgroup. Define
{\it the observable hull} $\hat H$ of $H$ as the minimal observable subgroup
of $G$ containing $H$. The stabilizer of any $H$-fixed vector in a
rational $G$-module contains $\hat H$. Therefore
$\kk[V]^H=\kk[V]^{\hat H}$ for any $G$-module $V$, and it is natural to
solve Hilbert's 14th problem for observable subgroups.

The following
famous theorem proved by F.~D.~Grosshans establishes a close connection
between Hilbert's 14th problem and the theory of affine embeddings.

\smallskip

\begin{theorem}~\cite{gr1}, \cite{gr}\label{gros}
Let $H$ be an observable subgroup of a reductive group $G$. The
following conditions are equivalent:

1) for any $G$-module $V$ the algebra $\kk[V]^H$ is finitely generated;

2) the algebra $\kk[G/H]$ is finitely generated;

3) there exists an affine embedding $G/H\hookrightarrow X$ such that
 $\codim_X (X\setminus (G/H))\ge 2$.
\end{theorem}

\begin{definition}
1)\  An observable subgroup $H$ in $G$ is said to be a \emph{Grosshans subgroup}
if $\kk[G/H]$ is finitely generated.

2)\ If $H$ is a Grosshans subgroup of $G$, then
$G/H\hookrightarrow X=\Sp\kk [G/H]$ is called
{\it the canonical embedding} of $G/H$, and $X$ is denoted by $CE(G/H)$.
\end{definition}

Note that any normal affine embedding $G/H\hookrightarrow X$ with
$\codim_X (X\setminus (G/H))\ge 2$ is $G$-isomorphic to the canonical
embedding~\cite{gr}. A homogeneous space $G/H$ admits such an embedding
if and only if $H$ is a Grosshans subgroup.

By Matsushima's criterion, $H$ is reductive if and only if $CE(G/H)=G/H$. For
non-reductive subgroups, $CE(G/H)$ is an interesting object canonically
associated with the pair $(G,H)$. It
allows to reformulate algebraic problems concerning the algebra $\kk[G/H]$ in
geometric terms.


\subsection{Popov-Pommerening's conjecture and Knop's theorem}
\label{subsec2.2}

\smallskip

\begin{theorem}~\cite{gr2}, \cite{don}, \cite[Th.16.4]{gr} \label{ttot}
Let $P^u$ be the unipotent radical of a parabolic subgroup $P$ of $G$.
Then $P^u$ is a Grosshans subgroup of $G$.
\end{theorem}

{\bf Proof.} \ Let $P=LP^u$ be a Levi decomposition and $U_1$ a maximal
unipotent subgroup of $L$. Then $U=U_1P^u$ is a maximal unipotent subgroup
of $G$, and $\kk[G]^U=(\kk[G]^{P^u})^{U_1}$. We know that $\kk[G]^U$ is finitely
generated (Theorem~\ref{teouinv}). On the other hand, Theorem~\ref{teouinv}
implies that
the $L$-algebra $\kk[G]^{P^u}$ is finitely generated if and only if
$(\kk[G]^{P^u})^{U_1}$ is, hence $\kk[G]^{P^u}$ is finitely generated.
(Another proof, using an explicit codimension 2 embedding,
is given in~\cite{gr2}.) $\diamondsuit$.

\smallskip

Let us say that a subgroup of a reductive group $G$
is {\it regular} if it is normalized
by a maximal torus in $G$. Generalizing Theorem~\ref{ttot},
V.~L.~Popov and K.~Pommerening conjectured that any observable regular
subgroup is a Grosshans subgroup.
At the moment,
a positive answer is known for groups $G$ of small rank~\cite{lt1},
\cite{lt2}, \cite{lt3}, and for some
special classes of regular subgroups (for example, for unipotent
radicals of parabolic subgroups of Levi subgroups of $G$~\cite{gr}).
Lin Tan~\cite{lt3} constructed explicitly canonical embeddings for
regular unipotent subgroups in $SL(n)$, $n\le 5$. A strong argument in favour
of Popov-Pommerening's conjecture is given in~\cite[Th.4.3]{bbk} in terms
of finite generation of induced modules, see also~\cite[\S\,23]{gr}.

\smallskip

Another powerful
method to check that the algebra $\kk[G/H]$ is finitely generated is provided by
the following theorem proved by F.~Knop.

\smallskip

\begin{theorem}~\cite{knH}, \cite{gr}
Suppose that $G$ acts on an irreducible
normal unirational variety $X$. If $c(X)\le 1$,
then the algebra $\kk[X]$ is finitely generated.
\end{theorem}

\begin{corollary}
If $H$ is observable in $G$ and $c(G/H)\le 1$, then $H$ is a Grosshans
subgroup.
\end{corollary}


\subsection{The canonical embedding of $G/P^u$}

Since the unipotent radical $P^u$ of a parabolic subgroup $P$
is a Grosshans subgroup of $G$, there exists a canonical embedding
$G/P^u\hookrightarrow CE(G/P^u)$. Such embeddings provide an interesting
class of affine factorial $G$-varieties, which was studied in~\cite{at2}.
Let us note that the Levi subgroup $L\subset P$ normalizes $P^u$, hence
acts $G$-equivariantly on $G/P^u$ and on $CE(G/P^u)$.
By $V_L(\lambda)$ denote a simple $L$-module with the highest weight $\lambda$.
Our approach
is based on the analysis of the $(G\times L)$-module decomposition of the
algebra $\kk[G/P^u]$ given by
$$
  \kk[G/P^u]=\bigoplus_{\lambda\in\Xi_+(G)} \kk[G/P^u]_{\lambda},
$$
where $\kk[G/P^u]_{\lambda}\cong V(\lambda)^*\otimes V_L(\lambda)$ is the
linear span of the matrix entries of the linear maps $V(\lambda)^{P^u}\to
V(\lambda)$ induced by $g\in G$, considered as regular functions on $G/P^u$.
(In fact, our method works for any affine embedding $G/P^u\hookrightarrow X$,
where $L$ acts $G$-equivariantly.) The multiplication structure looks like
$$
\kk[G/P^u]_{\lambda}\cdot\kk[G/P^u]_{\mu}=\kk[G/P^u]_{\lambda+\mu}\oplus
\bigoplus_i\kk[G/P^u]_{\lambda+\mu-\beta_i},
$$
where $\lambda+\mu-\beta_i$ runs over the highest weights of all "lower"
irreducible components in the $L$-module decomposition
$V_L(\lambda)\otimes V_L(\mu)=V_L(\lambda+\mu)\oplus\dots$.

\smallskip

Below we list the results from~\cite{at2}.

\smallskip

$\bullet$ Affine $(G\times L)$-embeddings $G/P^u\hookrightarrow X$ are
classified by finitely generated subsemigroups $S$ of $\Xi_+(G)$
having the property that all highest weights of the tensor product of
simple $L$-modules with highest weights in $S$ belong to $S$, too. Furthermore,
every choice of the generators $\lambda_1,\dots,\lambda_m\in S$ gives rise
to a natural $G$-equivariant embedding $X\hookrightarrow\Hom(V^{P^u}, V)$,
where $V$ is the sum of simple $G$-modules of highest weights $\lambda_1,\dots,
\lambda_m$. The convex cone $\Sigma^+$ spanned by $S$ is nothing else but the
dominant part of the cone $\Sigma$ spanned by the weight polytope of $V^{P^u}$.
In the case $X=CE(G/P^u)$, the semigroup $S$ coincides with $\Xi_+(G)$ and
$\Sigma$ is the span of the dominant Weyl chamber by the Weyl group of $L$.
In particular, if $G$ is
simply connected semisimple, then there is a natural inclusion
$$
  CE(G/P^u)\subset\bigoplus_{i=1}^l\Hom(V(\omega_i)^{P^u}, V(\omega_i)),
$$
where $\omega_1,\dots,\omega_l$ are the fundamental weights of $G$.

\smallskip

$\bullet$ The $(G\times L)$-orbits in $X$ are in bijection with the faces of
$\Sigma$ whose interiors contain dominant weights, the orbit representatives
being given by the projectors onto the subspaces of $V^{P^u}$ spanned
by eigenvectors of eigenweights in a given face. For the canonical embedding,
the $(G\times L)$-orbits correspond to the subdiagrams in the Dynkin diagram
of $G$ such that no connected component of such a subdiagram is contained in
the Dynkin diagram of $L$. Also we compute the stabilizers of points in
$G\times L$ and in $G$, and the modality of the action $G:X$.

\smallskip

$\bullet$ We classify smooth affine $(G\times L)$-embeddings
$G/P^u\hookrightarrow X$. In particular, the only non-trivial smooth canonical
embedding corresponds to $G=SL(n)$, $P$ is the stabilizer of a hyperplane
in $\kk^n$, and $CE(G/P^u)={\rm Mat}(n,n-1)$ with the $G$-action by left
multiplication.

\smallskip

$\bullet$ The techniques used in the description of affine $(G\times
L)$-embeddings of $G/P^u$ are parallel to those
developed in~\cite{dti} for the study of equivariant
compactifications of reductive groups. An analogy with monoids becomes
more transparent in view of the bijection between our affine
embeddings $G/P^u\hookrightarrow X$ and a class of algebraic monoids $M$
with the group of invertibles~$L$, given by $X=\Sp\kk[G\times^PM]$.

\smallskip

$\bullet$ Finally, we describe the $G$-module structure of the tangent
space of $CE(G/P^u)$ at the $G$-fixed point, assuming that
$G$ is simply connected simple. This space is obtained from
$\bigoplus_i\Hom(V(\omega_i)^{P^u}, V(\omega_i))$ by removing certain summands
according to an explicit algorithm.
The tangent space at the fixed point is at the same time the
minimal ambient $G$-module for $CE(G/P^u)$.


\subsection{Counterexamples}
The famous Nagata's counterexample to Hilbert's
14th problem~\cite{na} yields a 13-dimensional unipotent subgroup $H$
in $SL(32)$ acting naturally in $V=\kk^{32}$
such that the algebra of invariants $\kk[V]^H$ is not
finitely generated. This shows that the algebra $\kk[SL(32)/H]$ is not finitely
generated, or, equivalently, the complement to the open orbit in
any affine embedding $SL(32)/H\hookrightarrow X$ contains
a divisor.

Nagata's construction was simplified by R.~Steinberg. He proved
that $\kk[V]^H$ is not finitely generated for
the following 6-dimensional commutative unipotent linear group:

\smallskip

$$
  H=\left\{\left(\begin{array}{ccc}
                   \mbox{$\begin{array}{cc} 1 & c_1 \\ 0 & 1 \end{array}
                         $} &   & 0 \\
                    & \mbox{$\begin{array}{cc} . & . \\ . & . \end{array}
                          $} &  \\
           0 &  &  \mbox{$\begin{array}{cc} 1 & c_9 \\ 0 & 1 \end{array}
                         $} \\
                     \end{array} \right), \
                      \sum_{j=1}^9 a_{ij}c_j=0, \ i=1,2,3 \ \   \right\},
$$

\smallskip

\noindent
where the nine points $P_j=(a_{1j}:a_{2j}:a_{3j})$ are nonsingular points
on an irreducible cubic curve in the projective plane, their sum has infinite
order in the group of the curve, and $V=\kk^{18}$ (see~\cite{st} for details).

Another method to obtain counterexamples was proposed by P.~Roberts~\cite{ro}.
Consider the polynomial algebra $R=\kk[x,y,z,s,t,u,v]$ in 7 variables over
a not necessarily algebraically closed field $\kk$ of characteristic zero
with the grading $R=\oplus_{n\ge 0} R_n$
determined by assigning the degree $0$ to $x,y,z$ and the degree $1$ to
$s,t,u,v$. The elements
$s,t,u,v$ generate a free $R_0$-submodule in $R$ considered as $R_0$-module .
Choosing a natural number $m\ge 2$,
Roberts defines an $R_0$-module homomorphism on this submodule
$$
f: R_0s\oplus R_0t\oplus R_0u\oplus R_0v \to R_0
$$
given by $f(s)=x^{m+1}$, $f(t)=y^{m+1}$, $f(u)=z^{m+1}$, $f(v)=(xyz)^m$.
The submodule $\Ker\,f$ generates a subalgebra of $R$, which is denoted by $A$.
It is proved in~\cite{ro} that the $\kk$-algebra $B=R\cap QA$ is not finitely
generated. (Roberts shows how to construct an element in $B$ of any given degree
which is not in the subalgebra generated by elements of lower degree.)
A linear action of a 12-dimensional commutative unipotent group
on 19-dimensional vector space with the algebra of invariants
isomorphic to the polynomial algebra in one variable over $B$
is constructed in~\cite{ac}.

For a recent development in this direction, see~\cite{df}, \cite{f}.


\section{Some properties of affine embeddings}\label{Sec3}

\subsection{Affinely closed spaces and Luna's theorem}\label{subsecLu}

\begin{definition}\label{ac}
An affine homogeneous space $G/H$ of an affine algebraic group $G$ is called
\emph{affinely closed} if it admits only trivial affine embedding
$X=G/H$.
\end{definition}

Assume that $G$ is reductive.
Then $G/H$ is affinely closed implies $H$ is reductive.
By $N_G(H)$ (resp. $C_G(H)$) denote the normalizer (resp. the centralizer)
of $H$ in $G$, and by $W(H)$ denote the quotient $N_G(H)/H$.
It is known that $N_G(H)^0=H^0C_G(H)^0$ and both $N_G(H)$
and $C_G(H)$ are reductive~\cite[Lemma 1.1]{lr}.

\smallskip

\begin{theorem}~\cite{lu1}\label{afcl}
Let $H$ be a reductive subgroup of a reductive group $G$. The
homogeneous space $G/H$ is affinely closed if and only if the
group $W(H)$ is finite. Moreover, if $G$ acts
on an affine variety $X$ and the stabilizer of a point $x\in X$ contains a
reductive subgroup $H$ such that $W(H)$ is finite, then
the orbit $Gx$ is closed in $X$.
\end{theorem}

\begin{remark}
The last statement may be reformulated: if $H$ is reductive, the group
$W(H)$ is finite, and $H\subset H'\subset G$, where $H'$ is observable,
then $H'$ is reductive and $G/H'$ is affinely closed.
\end{remark}

\begin{remark}
Let $H$ be a Grosshans subgroup of $G$.
The following conditions are equivalent:

1)\ $H$ is reductive and $W(H)$ is finite;

2)\ $H$ is reductive and
for any one-parameter subgroup $\mu:\kk^*\to C_G(H)$ one has
$\mu(\kk^*)\subseteq H$;

3)\ the algebra $\kk[G/H]$ does not have non-trivial $G$-invariant ideals and
does not admit non-trivial $G$-invariant $\ZZ$-gradings;

4)\ the algebra $\kk[G/H]$ does not have non-trivial $G$-invariant ideals and
the group of $G$-equivariant automorphisms of $\kk[G/H]$ is finite.

5)\ no invariant subalgebra in $\kk[G/H]$ admits a non-trivial
$G$-invariant ideal.
\end{remark}

\begin{example}
1)\ Let $\rho:H\to SL(V)$ be an irreducible representation of a reductive
group $H$. Then the space $SL(V)/\rho(H)$ is affinely closed
($W(\rho(H))$ is finite by the Schur Lemma).

2)\ If $T$ is a maximal torus of $G$, then $W(T)$ is the Weyl group and
$G/T$ is affinely closed.
\end{example}

\begin{proposition}\label{prpol} Let $G$ be an affine algebraic group.
The following conditions are equivalent:

1) any monoid $S$ with $G(S)=G$ and $\overline{G(S)}=S$ coincides with $G$;

2) the group $G/G^u$ is semisimple.
\end{proposition}

{\bf Proof.}\ Let $G$ be reductive. The space $(G\times G)/\Delta(G)$
is affinely closed if and only if
the group $N_{G\times G}(\Delta(G))/\Delta(G)$ is finite.
But this is exactly the case when
the center of $G$ is finite. The same arguments work for any $G$
(Theorem~\ref{tm}). \ $\diamondsuit$

Below we give a proof of Theorem~\ref{afcl}
in terms of so-called adapted (or optimal)
one-parameter subgroups following G.~Kempf~\cite[Cor.4.5]{ke}.

We have to prove that if $G/H'$ is a quasi-affine
homogeneous space that is not affinely closed and $H\subset H'$ is a reductive
subgroup, then there  exists a one-parameter subgroup
$\nu: \kk^*\to C_G(H)$ such that $\nu(\kk^*)$ is not contained in $H$.
There is an affine
embedding $G/H'\hookrightarrow X$ with a $G$-fixed point $o$,
see~\ref{ssil}.
Denote by $x$ the image of $eH'$ in the open orbit on $X$.
By the Hilbert-Mumford criterion,
there exists a one-parameter subgroup $\gamma: \kk^*\to G$ such that
$\lim_{t\to 0} \gamma(t)x=o$. Moreover, there is a subgroup $\gamma$
that moves $x$ "most rapidly" toward $o$. Such $\gamma$ is called
{\it adapted to $x$},
for the precise definition see~\cite{ke}, \cite{pv2}.
For adapted $\gamma$, consider the parabolic subgroup
$$
  P(\gamma)=\{g\in G \mid \lim_{t\to 0} \gamma(t)g\gamma(t)^{-1}
  \ {\rm exists \ in} \ G \ \}.
$$
Then $P(\gamma)=L(\gamma)U(\gamma)$, where $L(\gamma)$ is the Levi
subgroup that is the centralizer of $\gamma(\kk^*)$ in $G$,
and $U(\gamma)$ is the unipotent radical of $P(\gamma)$.
By~\cite{ke}, \cite[Th.5.5]{pv2},
the stabilizer $G_x=H'$ is contained in
$P(\gamma)$. Hence there is an element $u\in U(\gamma)$ such that
$uHu^{-1}\subset L(\gamma)$.

We claim that $\gamma(\kk^*)$ is not
contained in $uHu^{-1}$. In fact, $\gamma$ is adapted to the element $ux$,
too~\cite[Th.3.4]{ke}, hence $\gamma(\kk^*)$ is not contained in the stabilizer of $ux$.
Thus $u^{-1}\gamma u$ is the desired subgroup $\nu$.

Conversely, suppose that there exists $\nu:\kk^*\to C_G(H)$ and $\nu(\kk^*)$
is not contained in $H$.
Consider the subgroup $H_1=\nu(\kk^*)H$.
The homogeneous fiber space $G*_{H_1}\kk$, where $H$ acts
on $\kk$ trivially and $H_1/H$ acts on $\kk$
by dilation, is a two-orbit embedding of $G/H$. $\diamondsuit$


\subsection{Affinely closed spaces in arbitrary characteristic}
In this subsection we assume that $\kk$ is an arbitrary algebraically
closed field. Suppose that $G$ acts on an affine
variety $X$. In positive characteristic, the structure of algebraic variety
on the orbit $Gx$ of a point $x\in X$ is not determined
(up to $G$-isomorphism) by the stabilizer $H=G_x$, and it is
natural to consider the isotropy subscheme $\tilde H$ at $x$, with $H$
as the reduced part, identifying $Gx$ and $G/\tilde H$. There is a
natural bijective purely inseparable and finite morphism $\pi:
G/H\to G/\tilde H$~\cite[4.3, 4.6]{hum}.
The following technical proposition shows
that this difficulty does not play an essential role for affinely closed spaces.

\smallskip

\begin{proposition}~\cite[Prop.8]{ar5}
The homogeneous space $G/H$ is affinely closed if and
only if $G/\tilde H$ is affinely closed.
\end{proposition}

\begin{definition}
We say that an affinely closed homogeneous space $G/H$ is
{\it strongly affinely
closed} if for any affine $G$-variety $X$
and any point $x\in X$ fixed by $H$ the orbit $Gx$ is closed in $X$.
\end{definition}

By Theorem~\ref{afcl}, in characteristic zero any affinely closed space is
strongly affinely closed.

The following notion was introduced by J.-P.Serre, cf.~\cite{ls}.

\begin{definition}
A subgroup $D\subset G$ is called
{\it $G$-completely reducible} ($G$-cr for short) if, whenever $D$ is
contained in a parabolic subgroup $P$ of $G$, it is contained in a
Levi subgroup of $P$.
\end{definition}

A $G$-cr subgroup is reductive.
For $G=GL(V)$ this notion agrees with the usual notion of complete reducibility.
In fact, if $G$ is any of the classical groups then the notions coincide,
although for symplectic and orthogonal groups this requires the assumption
that $\cchar \kk$ is a good prime for $G$.
The class of $G$-cr subgroups is wide. Some conditions which
guarantee that certain subgroups satisfy the $G$-cr condition can
be found in~\cite{mcn}, \cite{ls}.

The proof of Theorem~\ref{afcl} given above implies:

\smallskip

$\bullet$\ if $H$ is not contained in a proper parabolic subgroup of $G$,
then $G/H$ is strongly affinely closed;

\smallskip

$\bullet$\ if there exists $\nu:\kk^*\to C_G(H)$ such that $\nu(\kk^*)$
is not contained in $H$, then $G/H$ is not affinely closed;

\smallskip

$\bullet$\ if $H$ is a $G$-cr subgroup of $G$, then the following conditions
are equivalent:

1)\ $G/H$ is affinely closed;

2)\ $G/H$ is strongly affinely closed;

3)\ for any one-parameter subgroup $\nu:\kk^*\to C_G(H)$ one has
$\nu(\kk^*)\subseteq H$.

\begin{example}
The following example produced by George J. McNinch
shows that the group $W(H)$ may be
unipotent even for reductive $H$. Let $L$ be the space of
$(n\times n)$-matrices and $H$ the image of $SL(n)$ in
$G=SL(L)$ acting on $L$ by conjugation.

If $p=\cchar\kk \mid n$, then $L$ is an indecomposable $SL(n)$-module
with three composition factors, cf.~\cite[Prop.~4.6.10, a)]{mcn}. It
turns out that $C_G(H)^0$ is a one-dimensional unipotent group
consisting of operators of the form ${\mbox Id}+aT$, where $a\in \kk$,
and $T$ is a nilpotent operator on $L$ defined by $T(X)=\tr(X)E$.
The subgroup $H$ is contained in a quasi-parabolic
subgroup of $G$, hence $G/H$ is not strongly affinely closed.

In the simplest case $n=p=2$, we have $H\cong
PSL(2)\subset SL(4)$, $N_G(H)=HC_G(H)$ (because $H$ does not have
outer automorphisms), $C_G(H)$ is connected, and $W(H)\cong (\kk,+)$.
\end{example}

\smallskip

It would be very interesting to obtain a complete description of
affinely closed spaces in arbitrary characteristic and to answer
the following question: is it true that any affinely closed space is strongly
affinely closed ?


\subsection{Affinely closed spaces of non-reductive groups}

For non-reductive $G$, the class of affinely closed homogeneous spaces
is much wider. For example, it is well-known that an orbit of a unipotent
group acting on an affine variety is closed, hence any homogeneous space
of a unipotent group is affinely closed. Conversely, if any (quasiaffine)
homogeneous space of an affine group $G$ is affinely closed, then the connected
component of the identity in $G$ is
unipotent~\cite[10.1]{bir}, \cite[Th.4.2]{fs}.
In this subsection we give a
complete characterization of affinely closed homogeneous spaces of
non-reductive groups.

Let us fix the Levi decomposition $G=LG^u$ of the group $G$ in the
semidirect product of a reductive subgroup $L$ and the unipotent
radical $G^u$. By $\phi$ denote the homomorphism $G \to
G/G^u$. We shall identify the image of $\phi$ with $L$. Put
$K=\phi(H)$.

\begin{theorem} \cite[Th.2]{ate} \label{tm}
The following conditions are equivalent:

\ {\rm (1)}\ $G/H$ is affinely closed;

\ {\rm (2)}\ $L/K$ is affinely closed.
\end{theorem}

{\bf Proof.}\
The subgroup $H$ is observable in $G$
if and only if the subgroup $K$ is observable in $L$~\cite{su},
\cite[Th.7.3]{gr}.

\smallskip

Suppose that $L/K$ admits a non-trivial affine embedding. Then
there are an $L$-module $V$ and a vector $v\in V$ such that the
stabilizer $L_v$ equals $K$ and the orbit boundary $Y=Z\setminus
Lv$, where $Z=\overline{Lv}$, is nonempty. Let $I(Y)$ be the ideal
in $\kk[Z]$ defining the subvariety $Y$.
There exists a finite-dimensional $L$-submodule $V_1\subset
I(Y)$ that generates $I(Y)$ as an ideal. The inclusion
$V_1\subset\kk[Z]$ defines $L$-equivariant morphism $\psi:Z\to
V_1^*$ and $\psi^{-1}(0)=Y$. Then $L$-equivariant morphism $\xi:
Z\to V_2=V_1^*\oplus (V\otimes V_1^*)$, $z\to (\psi(z),
z\otimes\psi(z))$ maps $Y$ to the origin and is injective on the
open orbit in $Z$. Hence we obtain an embedding of $L/K$ in an
$L$-module such that the closure of the image of this embedding
contains the origin. Put $v_2=\xi(v)$. By the Hilbert-Mumford
Criterion, there is a one-parameter subgroup $\lambda: \kk^* \to
L$ such that $\lim_{t\to 0} \lambda(t)v_2=0$. Consider the weight
decomposition $v_2=v_2^{(i_1)}+\dots+v_2^{(i_s)}$ of the vector
$v_2$, where $\lambda(t)v_2^{(i_k)}=t^{i_k}v_2^{(i_k)}$. Here all
$i_k$ are positive.

By the identification $G/G^u=L$, one may consider $V_2$ as a
$G$-module. Let $W$ be a finite-dimensional $G$-module with a
vector $w$ whose stabilizer equals $H$. Replacing the pair $(W,
w)$ by the pair $(W\oplus (W\otimes W), w+w\otimes w)$, one may
suppose that the orbit $Gw$ intersects the line $\kk w$ only at
the point $w$. The weight decomposition shows that, for a sufficiently
large $N$, in the $G$-module
$W\otimes V_2^{\otimes N}$ one has $\lim_{t\to 0} \lambda(t)
(w\otimes v_2^{\otimes N})=0$ ($\lambda(\kk^*)$ may be considered
as a subgroup of $G$). On the other hand, the stabilizer of
$w\otimes v_2^{\otimes N}$ coincides with $H$. This implies that
the space $G/H$ is not affinely closed.

\smallskip

Conversely, suppose that $G/H$ admits a non-trivial affine
embedding. This embedding corresponds to a $G$-invariant
subalgebra $A\subset\kk[G/H]$ containing a non-trivial
$G$-invariant ideal $I$. Note that the algebra $\kk[L]$ may be
identified with the subalgebra in $\kk[G]$ of (left- or right-)
$G^u$-invariant functions, $\kk[G/H]$ is realized in $\kk[G]$ as
the subalgebra of right $H$-invariants, and $\kk[L/K]$ is the
subalgebra of left $G^u$-invariants in $\kk[G/H]$. Consider the
action of $G^u$ on the ideal $I$. By the Lie-Kolchin Theorem,
there is a non-zero $G^u$-invariant element in $I$. Thus the
subalgebra $A\cap\kk[L/K]$ contains the non-trivial $L$-invariant
ideal $I\cap\kk[L/K]$. If the space $L/K$ is affinely closed then
we get a contradiction with the following lemma.

\begin{lemma}\label{lv}
Let $L/K$ be an affinely closed space of a reductive group
$L$. Then any $L$-invariant subalgebra in $\kk[L/K]$ is finitely generated and
does not contain non-trivial $L$-invariant ideals.
\end{lemma}

{\bf Proof.}\
Let $B\subset\kk[L/K]$ be a non-finitely generated
invariant subalgebra. For any chain $W_1\subset W_2\subset
W_3\subset\dots$ of finite-dimensional $L$-invariant submodules in
$\kk[L/K]$ with $\cup_{i=1}^{\infty} W_i=\kk[L/K]$, the chain of
subalgebras $B_1\subset B_2\subset B_3\subset\dots$ generated by
$W_i$ does not stabilize. Hence one may suppose that all
inclusions here are strict. Let $Z_i$ be the affine $L$-variety
corresponding the algebra $B_i$. The inclusion $B_i\subset
\kk[L/K]$ induces the dominant morphism $L/K\to Z_i$ and
Theorem~\ref{afcl} implies that $Z_i=L/K_i$, $K\subset K_i$. But
$B_1\subset B_2\subset B_3\subset\dots$, and any $K_i$ is strictly
contained in $K_{i-1}$, a contradiction. This shows that $B$ is
finitely generated and, as proved above, $L$ acts transitively on
the affine variety $Z$ corresponding to $B$. But any non-trivial
$L$-invariant ideal in $B$ corresponds to a proper $L$-invariant
subvariety in $Z$.
\ $\diamondsuit$ \ \ Theorem~\ref{tm} is proved.
\ $\diamondsuit$

\begin{corollary}\label{c1}
Let $G/H$ be an affinely closed homogeneous space. Then for any affine
$G$-variety $X$ and a point $x\in X$ such that $Hx=x$, the orbit $Gx$
is closed.
\end{corollary}

{\bf Proof.}\
The stabilizer $G_x$ is observable in $G$, hence $\phi(G_x)$ is
observable in $L$. The subgroup $\phi(G_x)$ contains $K=\phi(H)$,
and Theorems~\ref{afcl},~\ref{tm} imply that the space $L/\phi(G_x)$ is
affinely closed. By Theorem~\ref{tm}, the space $G/G_x$ is
affinely closed.
\ $\diamondsuit$

\begin{corollary}\label{c2}
If $X$ is an affine $G$-variety and a point $x\in X$ is $T$-fixed,
where $T$ is a maximal torus of $G$, then the orbit $Gx$ is closed.
\end{corollary}

A characteristic-free description of affinely closed homogeneous spaces
for solvable groups is given in~\cite{te}.


\subsection{The Slice Theorem}\label{sslice}

The Slice Theorem due to D.~Luna~\cite{lu2} is one of the most important
technical tools in modern Invariant Theory. In this text we need only some
corollaries of the Slice Theorem related to affine embeddings~\cite{lu2},
\cite{pv2}.

\smallskip

$\bullet$\ Let $G/H\hookrightarrow X$ be an affine embedding with a closed
$G$-orbit isomorphic to $G/F$, where $F$ is reductive. By the Slice Theorem,
we may assume that $H\subseteq F$. Then there exists an affine embedding
$F/H\hookrightarrow Y$ with an $F$-fixed point such that $X$ is $G$-isomorphic
to the homogeneous fiber space $G*_F Y$. This allows to reduce many problems to
affine embeddings with a fixed point. On the other hand,
this gives us a $G$-equivariant
projection of $X$ onto $G/F$.

\smallskip

$\bullet$\ Let $G/H\hookrightarrow X$ be a smooth affine embedding with
closed $G$-orbit isomorphic to $G/F$. Then
$X$ is a homogeneous vector bundle over $G/F$. In particular,
if $X$ contains a $G$-fixed point, then $X$ is vector space with a linear
$G$-action.


\subsection{Fixed-point properties}\label{ssil}
Here we list some results concerning $G$-fixed points in affine embeddings.

\smallskip

$\bullet$\ If $G/H$ is a quasi-affine non affinely closed homogeneous
space, then $G/H$ admits an affine embedding with a $G$-fixed
point~\cite[Prop.3]{ar5}.

\smallskip

$\bullet$\  A homogeneous space $G/H$ admits an affine embedding
$G/H\hookrightarrow X$ such that $X=G/H\cup \{o\}$, where $o$ is
a $G$-fixed point, if and only if $H$ is a quasi-parabolic subgroup of
$G$~\cite[Th.4, Cor.5]{po}. In this case
the normalization of $X$ is an $HV$-variety and the normalization morphism
is bijective.

\smallskip

$\bullet$\ Consider the canonical decomposition
$\kk[G/H]=\kk\oplus \kk[G/H]_G$, where the first term corresponds to the
constant functions and $\kk[G/H]_G$ is the sum of all nontrivial
simple $G$-submodules in $\kk[G/H]$.
Suppose that $H$ is an observable subgroup of $G$.
The following conditions are equivalent~\cite[Prop.6]{ar5}:

(1)\ any affine embedding of $G/H$ contains a $G$-fixed point;

(2)\ $H$ is not contained in a proper reductive subgroup of $G$;

(3)\ $\kk[G/H]_G$ is an ideal in $\kk[G/H]$.

If $H$ is a Grosshans subgroup, then conditions (1)-(3) are equivalent to

(4)\ $CE(G/H)$ contains a $G$-fixed point.

\begin{example}
Let $G$ be a connected semisimple group and $P$ a parabolic subgroup
containing no simple components of $G$.
For $H=P^u$ the properties (1)-(4)
hold. In fact, (3) follows from the observation that $\kk[G/P^u]_G$ is the
positive part of a $G$-invariant grading
on $\kk[G/P^u]$ defined by the $G$-equivariant
action of a suitable one-parameter subgroup in the center of the Levi subgroup
of $P$ on $G/P^u$~\cite{ar5}.
\end{example}

\begin{proposition}
Let $H$ be an observable subgroup of $G$.

1) \ If either
$G/H$ is affinely closed or $H$ is a quasi-parabolic subgroup of $G$,
then $G/H$ admits only one normal affine embedding (up to
$G$-isomorphism);

2) \ if $G=\kk^*$ and $H$ is finite, then there exist only two normal
affine embeddings, namely $\kk^*/H$ and $\kk/H$;

3)\ in all other cases there exists an infinite
sequence
$$
 X_1\stackrel{\phi_1}{\longleftarrow}%
X_2\stackrel{\phi_2}{\longleftarrow} X_3%
\stackrel{\phi_3}{\longleftarrow}\dots
$$
of pairwise nonisomorphic normal
affine embeddings $G/H\hookrightarrow X_i$ and equivariant
dominant morphisms $\phi_i$.
\end{proposition}

{\bf Proof.}\
The statements are obvious for affinely closed $G/H$ and for $G=\kk^*$.
If $H$ is a quasi-parabolic subgroup, then $\kk[G/H]^U=\kk[t]$.
Suppose that $G/H\hookrightarrow X$ is a normal affine embedding.
Then $\kk[X]^U\subseteq\kk[t]$ is a graded integrally closed subalgebra
with $Q(\kk[X]^U)=\kk(t)$. This implies $\kk[X]^U=\kk[t]$ and
$\kk[X]=\kk[G/H]$, hence $X$ is $G$-isomorphic to the canonical embedding
of $G/H$.

In all other cases there exists an integrally
closed non-finitely generated invariant subalgebra $\B$ in $\kk[G/H]$
with $Q\B=\kk(G/H)$, see Proposition~\ref{nfg}. Let $f_1,f_2,\dots,f_n,
f_{n+1},\dots$ be a set of generators of $\B$
such that $\kk(f_1,\dots,f_n)=\kk(G/H)$.
Define $\B_k$ as the integral closure of
$\kk[\langle Gf_1,\dots,Gf_{n+k} \rangle]$ in
$\B$. The varieties $X_k=\Sp \B_k$ are birationally isomorphic to $G/H$ and
hence $G/H\hookrightarrow X_k$. Infinitely many of the $X_k$ are pairwise
nonisomorphic. Renumbering, one may suppose that all $X_k$ are nonisomorphic.
The chain
$$
  \B_1\subset\B_2\subset\B_3\dots
$$
corresponds to the desired chain
$$
 X_1\gets X_2\gets X_3\gets \dots \ \ \ \ \ \ \ \diamondsuit
$$


\section{Embeddings with a finite number of orbits}\label{Sec4}

\subsection{The characterization theorem} \label{subsub}
Spherical homogeneous spaces admit the following nice characterization
in terms of equivariant embeddings.

\smallskip

\begin{theorem}~\cite{ser},~\cite{lv},~\cite{akh1} \label{FO}
A homogeneous space $G/H$ is spherical if and only if
any embedding of $G/H$ has finitely many $G$-orbits.
\end{theorem}

To be more precise, F.~J.~Servedio proved that any affine spherical
variety contains finitely many $G$-orbits, D.~Luna, Th.~Vust
and D.~N.~Akhiezer extended this result to an arbitrary
spherical variety, and D.~N.~Akhiezer constructed
a projective embedding with infinitely many $G$-orbits
for any homogeneous space of positive complexity.

Now we are concerned with the following problem:
characterize all quasi-affine homogeneous spaces $G/H$ of a
reductive group $G$ with the property

\begin{center}
(AF) \quad
\itshape For any affine embedding $G/H\hookrightarrow X$, the
number of $G$-orbits in $X$ is finite.
\end{center}

\noindent It follows from results considered above that

1) spherical homogeneous spaces;

2) affinely closed homogeneous spaces;

3) homogeneous spaces of the group $SL(2)$

\noindent have property (AF).
Our main result in some sense gives a unification
of these three classes.

\begin{theorem}\cite{at}\label{main}
For a reductive subgroup $H\subseteq G$, (AF) holds if and only
if either $W(H)=N_G(H)/H$ is finite or any extension of $H$ by a
one-dimensional torus in $N_G(H)$ is spherical in~$G$.
\end{theorem}

\begin{corollary}\label{c>1}
For an affine homogeneous space $G/H$ of complexity $>1$, (AF)
holds if and only if $G/H$ is affinely closed.
\end{corollary}

\begin{corollary}\label{c=1}
An affine homogeneous space $G/H$ of complexity~$1$ satisfies (AF)
if and only if either $W(H)$ is finite, or $\rk W(H)=1$ and $N_G(H)$
is spherical.
\end{corollary}

\begin{corollary}\label{center}
Let $G$ be a reductive group with infinite center $Z(G)$ and
$H$ a reductive subgroup in $G$ that does not contain
$Z(G)^0$. Then  property (AF) holds for $G/H$ if and only if
$H$ is a spherical subgroup of $G$.
\end{corollary}

The proof of Theorem~\ref{main} is based on the analysis of
Akhiezer's construction \cite{akh1} of projective embeddings and on
some results of F.~Knop. We give this proof in~\ref{moda} obtaining
a more general Theorem~\ref{gener}.

Our method applied to an
arbitrary quasi-affine space $G/H$ gives a necessary condition
for property (AF) (see Remark~\ref{obs} below),
but a characterization of quasi-affine
spaces with property (AF) is not obtained yet.
Another open problem is to characterize Grosshans subgroups
$H$ of a reductive group $G$ such that $CE(G/H)$ contains only a finite
number of $G$-orbits~\cite{ar5}.


\subsection{Modality}\label{moda}

The aim of this subsection is to generalize Theorem~\ref{main} following
the ideas of~\cite{akh2}, and to find the maximal number of parameters
in a continuous family of $G$-orbits over all affine embeddings of a given
affine space $G/H$.

\begin{definition}\label{modality}
Let $F:X$ be an algebraic group action. The integer
$$
d_F(X)=\min_{x\in X}\codim\nolimits_X Fx=\trdeg \ \kk(X)^F
$$
is called the {\it generic modality} of the action. This is the number of
parameters in the family of generic orbits. The {\it modality} \
of $F:X$ is the integer $\modd_F X=\max_{Y\subseteq X}\ d_F(Y)$, where $Y$
runs through $F$-stable irreducible subvarieties of $X$.
\end{definition}

An action of modality zero is nothing else but an action with a finite number
of orbits. Note that $c(X)=d_B(X)$.
E.~B.~Vinberg~\cite{vi1} proved that $\modd_B(X)=c(X)$
for any $G$-variety $X$.
This means that if we pass from $X$ to a $B$-stable irreducible subvariety
$Y\subset X$, then the number of parameters for generic
$B$-orbits does not increase.
Simple examples show that
the inequality $d_G(X)\le \modd_G(X)$ can be strict.
This motivates the following

\begin{definition}\label{assa}
With any $G$-variety $X$ we associate the integer
$$
m_G(X)=\max\nolimits_{X'}\modd_G(X'),
$$
where $X'$ runs through all $G$-varieties birationally $G$-isomorphic
to $X$.
\end{definition}

For a homogeneous space $G/H$ we have
$m_G(G/H)=\max_X\modd_G(X)$, where $X$ runs through
all embeddings of $G/H$.

It is clear that for any subgroup $F\subset G$ the inequality
$m_G(X)\le m_F(X)$ holds. In particular, $m_G(X)\le c(X)$.  The next
theorem shows that $m_G(X)=c(X)$.

\begin{theorem} \cite{akh2} \label{akhi2}
There exists a projective $G$-variety $X'$ birationally
$G$-isomorphic to $X$ such that $\modd_G(X')=c(X)$.
\end{theorem}

Now we introduce an affine counterpart of $m_G(X)$.

\begin{definition}\label{bssb}
With any quasi-affine homogeneous space $G/H$ we associate the integer
$$
a_G(G/H)=\max\nolimits_X \modd_G(X),
$$
where $X$ runs through all affine embeddings $G/H\hookrightarrow X$.
\end{definition}

\begin{theorem}\cite{ar3}\label{gener}
Let $H$ be a reductive subgroup of $G$.

(1) If the group $W(H)$ is finite, then $a_G(G/H)=0$;

(2) If $W(H)$ is infinite, then
$$
a_G(G/H)=\max\nolimits_{H_1}\ c(G/H_1),
$$
where $H_1$ runs through all non-trivial extensions of $H$ by a one-dimensional
subtorus of $C_G(H)$.
In particular, $a_G(G/H)=c(G/H)$ or $c(G/H)-1$.
\end{theorem}

{\bf Proof.}\ {\bf Step 1} -- Affine cones. Consider the natural
surjection $\kappa: N_G(H)\to W(H)$.

\begin{proposition}\label{basic}
Let $H$ be an observable subgroup of $G$.
Suppose that there is a non-trivial one-parameter subgroup
$\lambda: \kk^*\to W(H)$ and put
$H_1=\kappa^{-1}(\lambda(\kk^*))$. Then there exists an affine embedding
$G/H\hookrightarrow X$ with $\modd_G(X)\ge c(G/H_1)$.
\end{proposition}

The idea of the proof is to apply Akhiezer's construction~\cite{akh2}
to the homogeneous space $G/H_1$ and to consider the affine cone over a
projective embedding $G/H_1\hookrightarrow X'$ with $\modd_G(X')=c(G/H_1)$

\begin{lemma}\label{eigenvector}
In notation of Proposition~\ref{basic}, there exists a finite-dimensional
$G$-module $V$ and an $H_1$-eigenvector $v\in V$ such that

1) the orbit $G\langle v\rangle$ of the line $\langle v\rangle$
in $\PP(V)$ is isomorphic to $G/H_1$;

2) $H$ fixes $v$;

3) $H_1$ acts transitively on $\kk^*v$;

4) $\modd_G(\overline{G\langle v\rangle})=c(G/H_1)$.
\end{lemma}

{\bf Proof of Lemma~\ref{eigenvector}.}\ By Chevalley's theorem,
there exist a $G$-module $V'$ and a vector
$v'\in V'$ having property 1).
Let $\chi$ be the eigenweight of $H$ at~$v'$. Since
$H$ is observable in $G$, each finite-dimensional $H$-module
can be embedded into a finite-dimensional $G$-module~\cite{bbhm}.
In particular, there exists a $G$-module
$V''$ containing $H$-eigenvectors of the weight $-\chi$.
Among them we can choose an $H_1$-eigenvector
$v''$ and set
$V=V'\otimes V''$, $v=v'\otimes v''$. This pair has properties
1)-2).

If $H_1$ does not act transitively on $\kk^*v$, then take an
arbitrary $G$-module $W$ containing a vector with stabilizer $H$.
Take an $H_1$-eigenvector in $W^H$ with nonzero weight and replace $V$
by $V\otimes W$ and $v$ by $v\otimes w$. Conditions 1)-3) are now satisfied.

By a result of Akhiezer~\cite{akh2}, we can find a pair $(V',v')$
with properties 1) and 4).  Then we proceed as above obtaining a pair
$(V,v)$. The closure
$\overline{G\langle v\rangle}\subseteq\PP(V)$ lies in the image of the Segre
embedding

$$
\PP(V')\times\PP(V'')\times\PP(W)\hookrightarrow\PP(V),
$$
and it projects $G$-equivariantly onto
$\overline{G\langle v'\rangle}\subseteq\PP(V')$.
Now properties 1)-4) are satisfied for the pair $(V,v)$.
$\diamondsuit$

\begin{remark}\label{H=U}
If $H$ is reductive, then one can find $v$ in Lemma~\ref{eigenvector}
such that $G_v=H$. This is not
possible for an arbitrary observable subgroup, see~\cite[Remark~2]{at}.
\end{remark}

{\bf Proof of Proposition~\ref{basic}.}\
Let $(V,v)$ be the pair from Lemma~\ref{eigenvector}. Put $H'=G_v$
and  $\tilde X=\overline{Gv}$.
By properties 1)-3) and since $H_1/H$ is isomorphic to $\kk^*$,
$H'$ is a finite extension of $H$. By 3), the closure of the orbit $Gv$
in $V$ is a cone, therefore 4) implies the inequality
$\modd_G(\tilde X)\ge c(G/H_1)$.

Consider now the morphism $G/H\to G/H'$. It determines an embedding
$\kk[G/H']\subseteq \kk[G/H]$.
Let $A$ be the integral closure of the subalgebra
$\kk[\tilde X]\subseteq \kk[G/H']$ in the field $\kk(G/H)$.
We have the following commutative diagrams:
$$
\begin{array}{ccccccccc}
A & \hookrightarrow & \kk[G/H] & \hookrightarrow &  \kk(G/H) &
\ \ \ \ \ \ \ & \Sp A & \hookleftarrow & G/H \\
\uparrow &       & \uparrow &  & \uparrow &
\ \ \ \ \ \ \ & \downarrow & & \downarrow \\
\kk[\tilde X]       & \hookrightarrow & \kk[G/H'] & \hookrightarrow &
\kk(G/H') & \ \ \ \ \ \ \ & \tilde X & \hookleftarrow & G/H'
\end{array}
$$

The affine variety $X=\Sp A$ with the natural $G$-action can be regarded as
an affine embedding of
$G/H$. The embedding $\kk[\tilde X]\subseteq A$ defines a finite
(surjective) morphism $X \to \tilde X$, therefore
$\modd_G(X)=\modd_G(\tilde X)\ge c(G/H_1)$.
$\diamondsuit$

\medskip

{\bf Step 2.}\  Here we formulate several results due to F.~Knop.

\begin{lemma}[{\cite[7.3.1]{kn2}, see also~\cite[Lemma~3]{at}}]
\label{B-inv}
Let $X$ be an irreducible $G$-variety, and $v$ a $G$-invariant
valuation of $\kk(X)$ over $\kk$ with residue field~$\kk(v)$. Then
$\kk(v)^B$ is the residue field of the restriction of $v$ to~$\kk(X)^B$.
\end{lemma}

\smallskip

\begin{definition}~\cite[\S 7]{kn3}
\label{source}
Let $X$ be a normal $G$-variety. A discrete $\QQ$-valued
$G$-invariant valuation of $\kk(X)$ is said to be \emph{central}
if it vanishes on $\kk(X)^B\setminus\{0\}$. A \emph{source} of
$X$ is a non-empty $G$-stable subvariety $Y\subseteq X$ that is
the center of a central valuation of~$\kk(X)$.
\end{definition}

The following lemma is an easy consequence of~\cite{kn3},
for more details see~\cite[Lemma 4]{at}.

\begin{lemma}
\label{1-param}
If $X$ is a normal affine $G$-variety containing a proper source, then
there exists a one-dimensional torus
$S\subseteq\Aut_G(X)$ such that $\kk(X)^B\subseteq\kk(X)^S$.
(Here $\Aut_G(X)$ is the group of
$G$-equivariant automorphisms of $X$).
\end{lemma}

{\bf Step 3.}
Assertion (1) of Theorem~\ref{gener} follows from Theorem~\ref{afcl}.
To prove (2) we use Proposition~\ref{basic}. Since $H$ is reductive,
the group $W(H)$ is reductive and contains a one-dimensional subtorus
$\lambda(\kk^*)$. Hence $a_G(G/H)\ge c(G/H_1)\ge c(G/H)-1$.
If there exists a one-dimensional torus
in $W(H)$ such that $c(G/H)=c(G/H_1)$, we obtain an affine embedding of
$G/H$ of modality $c(G/H)$.

Conversely, let $G/H\hookrightarrow X$ be an affine embedding of modality
$c(G/H)$. We have to find a one-dimensional subtorus
$\lambda(\kk^*)\subseteq W(H)$ such that $c(G/H_1)=c(G/H)$.
By the definition of modality, there exists a proper
$G$-invariant subvariety $Y\subset X$ such that the codimension
of a generic $G$-orbit in $Y$ is $c(G/H)$, hence $c(Y)=c(G/H)$.
Consider a $G$-invariant valuation $v$ of $\kk(X)$ with center
$Y$. For the residue field $\kk(v)$ we have
$\trdeg \kk(v)^B\ge \trdeg \kk(Y)^B$,
therefore $\trdeg \kk(v)^B=\trdeg \kk(X)^B$. If the restriction of
$v$ to $\kk(X)^B$ is non-trivial, then, by Lemma~\ref{B-inv},
$\trdeg \kk(v)^B<\trdeg \kk(X)^B$, a contradiction. Thus, $v$
is central and
$Y$ is a source of $X$. Lemma~\ref{1-param} provides a one-dimensional
subtorus $S\subseteq Aut_G(X)\subseteq Aut_G(G/H)=W(H)$ that yields
an extension of $H$ of the same complexity.
$\diamondsuit$

\smallskip

Note that Theorem~\ref{main} is a particular case of Theorem~\ref{gener} with
$a_G(G/H)=0$.

\begin{remark}\label{obs}
If $H$ is an observable subgroup and $W(H)$ contains a non-trivial subtorus,
then the formula $a_G(G/H)=\max_{H_1} c(G/H_1)$ can be obtained by the same
arguments. In particular, Corollary~\ref{center} holds for observable $H$.
But for non-reductive $H$ the group $W(H)$ can be unipotent~\cite{at}:
this is the case when $G=SL(3)\times SL(3)$ and

\smallskip
$$
 H=\left\{ \left( \begin{array}{ccc}
                  1 & a & b+\frac{a^2}{2} \\
                  0 & 1 & a \\
                  0 & 0 & 1
                  \end{array}
  \right), \left( \begin{array}{ccc}
                  1 & b & a+\frac{b^2}{2} \\
                  0 & 1 & b \\
                  0 & 0 & 1
                  \end{array} \right) \left| \right. \ a,\ b \in \kk
                  \right\}.
$$
\smallskip
For such subgroups our proof yields only the inequality $a_G(G/H)\le c(G/H)-1$.
\end{remark}

Let us mention an application of Theorem~\ref{gener} which may be regarded
as its algebraic reformulation.
Let $G$ be a connected semisimple group.
Note that, for the action by left multiplication,
one has $c(G)=\frac{1}{2}\,(\dim\,G - \rk G)$ and
$c(G/S)=\frac{1}{2}\,(\dim\,G - \rk G)-1$, where $S$ is a one-dimensional
subtorus in $G$.
Applying Theorem~\ref{gener} to the case $H=\{ e\}$, we obtain

\smallskip

\begin{theorem}~\cite{ar4}
Let $A\subset\kk[G]$ be a left $G$-invariant finitely generated subalgebra
and $I\subset A$ a $G$-invariant prime ideal. Then
$$
\trdeg (Q(A/I))^G \ \le \ \frac{1}{2}\,(\dim\,G - \rk G)-1. \eqno(1)
$$
Moreover, there exist a subalgebra $A$ and an ideal $I$ such that $(1)$
is an equality.
\end{theorem}

\begin{example}
The closure of an $SL(3)$-orbit in an algebraic $SL(3)$-variety $X$ may contain
at most 3-parameter family of $SL(3)$-orbits. If $X$ is affine, then
the maximal number of parameters equals 2.
\end{example}


\subsection{Equivariant automorphisms and symmetric embeddings}
The group $\Aut_G(G/H)$ of $G$-equivariant automorphisms of
$G/H$ is isomorphic to $W(H)$. The action $W(H):G/H$
is induced by the action $N_G(H):G/H$ by right multiplication, i.e.,
$n*gH=gn^{-1}H$.
Let $G/H\hookrightarrow X$ be an embedding. The group
$\Aut_G X$ preserves the open orbit, and may be considered as a subgroup
of $W(H)$.

\begin{definition}\label{vs-def}
An embedding $G/H\hookrightarrow X$ is said to be
\emph{symmetric} if $W(H)^0\subseteq\Aut_G(X)$. If $\Aut_G(X)=W(H)$, we
say that $X$ is {\it very symmetric}.
\end{definition}

\begin{lemma} The following affine embeddings are very symmetric:

1) an affine embedding of a spherical homogeneous space;

2) the canonical embedding $CE(G/H)$;

3) an affine monoid $M$ considered as the embedding
$G(M)/\{e\}\hookrightarrow M$.
\end{lemma}

{\bf Proof.}\ 1) Let $G/H$ be a quasi-affine spherical homogeneous space.
By the Schur Lemma, the group $W(H)$ acts on any isotypic component
of $\kk[G/H]$ by dilation.
Hence any $G$-invariant subspace of $\kk[G/H]$ is also $W(H)$-invariant.

2) The group $W(H)$ acts on $G/H$ and on $\kk[G/H]$, thus on $\Sp\kk[G/H]$.

3) The group $W(H)\cong G(M)$ acts on $M$ by right multiplication.
$\diamondsuit$

\begin{proposition}\label{pprr}
Let $H$ be a reductive subgroup of $G$. The following conditions are equivalent:

(1) there exists a unique symmetric embedding $X=G/H$;

(2) $W(H)^0$ is a semisimple group.
\end{proposition}

{\bf Proof.}\ The existence of a non-trivial affine embedding
$G/H\hookrightarrow X$ with $\dim\Aut_G(X)=\dim W(H)$ means that
$G/H$ as a $(G\times W(H)^0)$-homogeneous space is not affinely closed.
By $L$ denote the $(G\times W(H)^0)$-stabilizer of the point  $eH$.
Then $L=\{(n, nH) \mid n\in\kappa^{-1}(W(H)^0) \}$ and
the group $N_{G\times W(H)^0}(L)/L$ is finite if and only if $W(H)$
is semisimple. $\diamondsuit$

\smallskip

Proposition~\ref{pprr} implies that
in the case of affine $SL(2)$-embeddings
only the trivial embedding $X=SL(2)$ is symmetric. In fact, in all other
cases with normal $X$ the group $\Aut_{SL(2)} X$ is a Borel
subgroup of $SL(2)$~\cite[III.4.8, Satz 1]{kr}.
The theorem below is a partial generalization of this result.

\smallskip

\begin{theorem}~\cite{at2} \label{teoavt}
Let $G/H\hookrightarrow X$ be an affine embedding with a finite number of
$G$-orbits and with a $G$-fixed point. Then the group $\Aut_G(X)^0$
is solvable.
\end{theorem}

We begin the proof with the following

\begin{lemma}\label{lem1}
Let $X$ be an affine variety with an action of a connected
semisimple group $S$. Suppose that there is a point $x\in X$
and a one-parameter subgroup $\gamma:\kk^*\to S$ such that
$\lim_{t\to 0}\delta(t)x$ exists in $X$ for any subgroup $\delta$
conjugate to~$\gamma$. Then $x$ is a $\gamma(\kk^*)$-fixed
point.
\end{lemma}

{\bf Proof.}\
Let $T$ be a maximal torus in $S$ containing
$\gamma(\kk^*)$. One can realize $X$ as a closed $S$-stable
subvariety in $V$ for a suitable $S$-module $V$.  Let
$x=x_{\lambda_1}+\dots+x_{\lambda_n}$ be the weight
decomposition (with respect to~$T$) of $x$ with weights
$\lambda_1,\dots,\lambda_n$.  One-parameter subgroups of $T$
form the lattice $\Xi_{*}(T)$ dual to the character lattice
$\Xi(T)$. The existence of $\lim_{t\to 0}\gamma(t)x$ in $X$ means
that all pairings $\langle\gamma,\lambda_i\rangle$ are
non-negative.  Let $\gamma_1,\dots,\gamma_m$ be all the
translates of $\gamma$ under the action of the Weyl group
$W=N_S(T)/T$. By assumption,
$\langle\gamma_j,\lambda_i\rangle\ge 0$ for any $i=1,\dots,n$,
$j=1,\dots,m$, hence
$\langle\gamma_1+\dots+\gamma_m,\lambda_i\rangle\ge 0$.  Since
$\gamma_1+\dots+\gamma_m=0$, one has
$\langle\gamma_j,\lambda_i\rangle=0$ for all $i,\ j$. This shows
that the points~$x_{\lambda_i}$ (and~$x$) are
$\gamma(\kk^*)$-fixed. \ $\diamondsuit$

\smallskip

The next
proposition is a generalization of~\cite[Th.4.3]{gr2}.

\begin{proposition} \label{prop1}
Suppose that $G/H\hookrightarrow X$ is an affine embedding with
a non-trivial $G$-equivariant action of a connected semisimple
group~$S$. Then the orbit $S*x$ is closed in $X$, $\forall x\in G/H$.
\end{proposition}

{\bf Proof.}\
We may assume $x=eH$. If $S*x$ is not closed, then,
by~\cite[Th.1.4]{ke}, there is a one-parameter subgroup
$\gamma:\kk^*\to S$ such that the limit
$$
  \lim_{t\to 0}\gamma(t)*x
$$
exists in $X$ and does not belong to $S*x$. Replacing $S$ by a
finite cover, we may assume that $S$ embeds in $N_G(H)$ (and
thus in~$G$) with a finite intersection with~$H$. By the
definition of $*$-action, one has
$\gamma(t)*x=\gamma(t^{-1})x$. For any $s\in S$ the limit
$$
\lim_{t\to 0} (s\gamma(t))*x=\lim_{t\to 0}\gamma(t^{-1})s^{-1}x
$$
exists. Hence
$\lim_{t\to 0}s\gamma(t^{-1})s^{-1}x$ exists, too.
This shows that
for any one-parameter subgroup $\delta$ of $S$,
conjugate to $-\gamma$,
$\lim_{t\to 0}\delta(t)x$ exists in~$X$.
Lemma~\ref{lem1} implies that $x=\lim_{t\to 0}\gamma(t)*x$,
and this contradiction proves Proposition~\ref{prop1}. \ $\diamondsuit$

\smallskip

{\bf Proof of Theorem~\ref{teoavt}.}\
Suppose that $\Aut_G(X)^0$ is not solvable. Then there is
a connected semisimple group $S$ acting on $X$ $G$-equivariantly.
By Proposition~\ref{prop1}, any $(S,*)$-orbit in the open
$G$-orbit of $X$ is closed in~$X$.

Let $X_1$ be the closure of a $G$-orbit in~$X$. Since
$G$ has a finite number of orbits in $X$, the variety $X_1$ is
$(S,*)$-stable. Applying the above arguments to~$X_1$, we
show that any $(S,*)$-orbit in $X$ is closed. But in this case
all $(S,*)$-orbits have the same dimension $\dim S$. On the
other hand, a $G$-fixed point is an $(S,*)$-orbit, a
contradiction. \ $\diamondsuit$

\begin{corollary}[of the proof]
Let $X$ be an affine $G$-variety with an open $G$-orbit. Suppose that

(1) \ a semisimple group $S$ acts on $X$ effectively and $G$-equivariantly;

(2) \ the dimension of a closed $G$-orbit in $X$ is less than $\dim S$.

Then the number of $G$-orbits in $X$ is infinite.
\end{corollary}

\begin{corollary}
Let $M$ be a reductive algebraic monoid with zero. Then the number of
left (right) $G(M)$-cosets in $M$ is finite if and only if $M$ is commutative.
\end{corollary}

The following corollary gives a partial answer to a question posed in
subsection~\ref{subsub}.

\begin{corollary}
The number of $G$-orbits in $CE(G/P^u)$ is finite if and only if
either $P\cap G_i=G_i$ or $P\cap G_i=B\cap G_i$ for each simple factor
$G_i\subseteq G$.
\end{corollary}

In many cases, Theorem~\ref{teoavt} may be used to show that the group
$\Aut_G(X)$ cannot be very big. On the other hand, the group
$\Aut_G(X)$ may be finite (trivial), in particular, for $X=G/H$ with
affinely closed $G/H$. Answering a question from~\cite{at}, I.~V.~Losev
proposed an example of an observable non-reductive subgroup $H$ in $SL(n)$,
where $W(H)$ is finite. (This example is included in the electronic
version of~\cite{at2}.) Note that any affine embedding of $SL(n)/H$ gives
an example of a locally transitive non-transitive reductive group action
on an affine variety with a finite group of equivariant automorphisms.
\smallskip

Finally, we give a variant of Theorem~\ref{main} for symmetric
embeddings.

\smallskip

\begin{theorem}~\cite[Prop.2]{at} \label{vs}
Let $H$ be a reductive subgroup of
$G$. Every symmetric affine embedding of $G/H$ has
finitely many $G$-orbits if and only if
either (AF) holds or $W(H)^0$ is semisimple.
\end{theorem}


\section{Application One: Invariant algebras on homogeneous spaces
of compact Lie groups}\label{Sec5}


\subsection{Invariant algebras and self-conjugate algebras}\label{ndss}

For any compact topological space
$M$ the set $C(M)$ of all continuous $\CC$-valued functions
on $M$ is a commutative Banach algebra with respect to pointwise addition,
multplication, and the uniform norm. We shall consider the case, where
$M=K/L$ is a
homogeneous space of a compact connected Lie group $K$.
Let us recall that $A$ is an {\it invariant algebra} on $M$
if $A$ is a $K$-invariant uniformly closed subalgebra with unit in $C(M)$.
In this section $G$ (resp. $H$) denotes the complexification of $K$
(resp. $L$). The group $G$ is a complex reductive algebraic group
with a reductive subgroup $H$.

The main problem is to describe all invariant algebras on a given space $M$
and to study their properties. Let us start with a particular class of
invariant algebras.

\begin{definition}
An invariant algebra $A$ is {\it self-conjugate} if $f\in A$ implies
$\overline f\in A$, where the bar denotes the complex conjugation.
\end{definition}

The classification of self-conjugate invariant algebras is based on the
Stone-Weierstrass Theorem. Here we follow~\cite{latdis}.

\medskip

{\bf The Stone-Weierstrass Theorem.} \ {\it Let $R$ be a compact
topological space and $A$
a subalgebra with unit in $C(R)$ such that

\ 1) \ $A$ separates points on $R$, i.e., for any $x_1\ne x_2\in R$ there exists
$f\in A$ such that $f(x_1)\ne f(x_2)$;

\ 2) \ $A$ is invariant under complex conjugation.

Then $A$ is dense in $C(R)$.}

\medskip

Given a self-conjugate invariant algebra $A$, define an equivalence
relation on $M$:
$x\sim y$ if and only if $f(x)=f(y)$ for any $f\in A$. The space $M'$ of
equivalence classes is a homogeneous $K$-space, hence $M'=K/L'$, where $L'$
is a closed subgroup containing $L$.
By construction, the self-conjugate algebra $A$ separates points on $M'$ and
thus $A=C(M')$. Conversely, for any $L\subseteq L'\subseteq K$ the inverse
image of $C(K/L')$ under the projection $K/L\to K/L'$ determines a self-conjugate
invariant algebra on $M$.
This shows that self-conjugate invariant algebras on $M$
are in one-to-one correspondence with closed subgroups $L'$, \
$L\subseteq L'\subseteq K$.


\subsection{Spherical functions}
The space $M=K/L$ may be considered as a compact subset of the affine
homogeneous space $X_0=G/H$. Moreover, $M$ is a real
form of $X_0$ in the natural sense. In particular, the restriction of
polynomial functions to $M$ determines an embedding
$\CC\,[X_0]\hookrightarrow C(M)$. Denote the image of this embedding
by $\CC\,[M]$.

\begin{definition}
A function $f\in C(M)$ is called {\it spherical} if the linear span
$\langle Kf \rangle$ is
finite-dimensional. More generally, for a linear action of a Lie group $K$
on vector space $V$,
a vector $v\in V$ is {\it spherical} if $\dim \langle Kv \rangle <\infty$.
\end{definition}

Denote by $V_{sph}$ the subspace of all spherical vectors in $V$.

\begin{proposition}
The algebra $\CC\,[M]$ coincides with $C(M)_{sph}$.
\end{proposition}

{\bf Proof.} Any regular function is contained in
finite-dimensional invariant subspace. Conversely,
any complex finite-dimensional representation of $K$ is completely
reducible and any irreducible component may be considered as a simple
$G$-module. Hence the matrix entries of such a module are in $\CC\,[M]$.
If $f\in C(M)$ is spherical and
$V=\langle Kf\rangle$, then $f$ is a linear combination
of the matrix entries of the dual representation $K:V^*$.
Indeed, let $f_1,\dots,f_k$ be a basis in $V$.
For any $f\in V, \ g\in K$ one has
$f_i(g^{-1}eL)=\sum a_{ij}(g)f_j(eL)$ and $f_i(gL)=\sum c_ja_{ij}(g^{-1})$,
where $c_j=f_j(eL)$ are constants. $\diamondsuit$

\smallskip

By the Peter-Weyl Theorem, the matrix entries
(with respect to some orthonormal basis) over all irreducible
finite-dimensional representations of $K$ form an orthonormal basis
in space $L^2(K)$.
Spherical functions are finite linear combinations of the basic elements. They
form uniformly dense subspace in $C(K)$. The following generalization
of this result plays a key role in this section.

\smallskip

\begin{proposition}~\cite[Th.5.1]{pw1},~\cite[2.16]{pw2}.
Given a continuous linear representation of a compact Lie group $K$
in Fr\'echet space $E$, the subspace $E_{sph}$ is dense in $E$.
\end{proposition}

In particular, in any invariant algebra spherical functions form a dense
subalgebra. Moreover, if $S$ is $K$-invariant subspace in $C(M)_{sph}$
and $\overline{S}$ is its uniform closure in $C(M)$, then
$\overline{S}\cap C_{sph}(M)=S$. (For the proof see~\cite[Lemma~14]{gl}.)
Finally, we get

\begin{theorem}\label{t1}
There is a natural bijection $\psi$
between invariant algebras on the space $M$
and invariant subalgebras in $\CC\,[M]$.
More precisely, $\psi(A)=\A=A_{sph}=A\cap\CC\,[M]$ and
$\psi^{-1}(\A)=\overline{\A}$.
\end{theorem}

This result provides nice connections between functional and algebraic
problems. To make this link really useful we need to reformulate functional
properties in algebraic terms and back. For this purpose we are going to use
geometric language of affine embeddings.


\subsection{Finitely generated invariant algebras and affine embeddings}

\begin{definition}
An invariant algebra $A$ is finitely generated if it is generated (as a
Banach algebra) by $K$-invariant finite-dimensional subspace.
\end{definition}

An invariant algebra $A$ is finitely generated if and only if $A_{sph}$
is a finitely
generated algebra. It is clear that $C(M)$ is finitely generated.
As follows from the discussion above, any self-conjugate invariant algebra
is finitely generated. The question when any invariant subalgebra in
$\CC\,[M]$ is finitely generated will be considered in the last section.

Any finitely generated subalgebra $\A\subset\CC\,[G/H]$ defines an affine
$G$-variety $X=\Sp \A$ with an open orbit isomorphic to $G/F$, where
$F$ is an observable subgroup containing $H$. The inclusion $\A\subset\CC\,[G/H]$
defines the morphism $\phi:G/H\to X$ and the base point $x_0=\phi(eH)$.
If we look at $\A$ as at an abstract $G$-algebra, then there may exist
different equivariant inclusion homomorphisms
$\A\to\CC\,[G/H]$ with the same image.
Two different base points $x_0\in X$ and $x_0'\in X$ determine the same
subalgebra $\A\subset \CC\,[G/H]$ if and only if
there exists $n\in\Aut_G(X)$ such that
$x_0=nx_0'$. (Corresponding inclusions $\A\subset\CC\,[G/H]$ differ by a
$G$-equivariant automorphism of $\A$.)
Let us denote the subalgebra $\A$ as $\A(X,x_0)$ and the
corresponding invariant algebra $\overline{\A(X,x_0)}$ as $A(X, x_0)$.
We have proved:

\begin{theorem}\label{t2}
Invariant finitely generated algebras on the space $M=K/L$ are in one-to-one
correspondence with the following data:

1) an affine embedding $G/F\hookrightarrow X$, where $F\subseteq G$
is an observable subgroup containing $H$;

2) an $H$-fixed point $x_0$
in the open $G$-orbit on $X$, which is defined
up to the action of $\Aut_G(X)$.
\end{theorem}

It is natural to classify invariant algebras up to some equivalence.
The group of $K$-equivariant automorphisms of $M$ is the group $N=N_K(L)/L$,
acting as $n*kL=kn^{-1}L$. This action defines a $K$-equivariant action
$N:C(M)$. The group $N$ acts transitively on the set $M^L$.

\begin{definition}
Two invariant algebras $A_1$ and $A_2$ on $M$ are {\it equivalent}
if there exists $n\in N$ such that $n*A_1=A_2$.
\end{definition}

Clearly, this equivalence preserves all reasonable properties of invariant
algebras. In terms of Theorem~\ref{t2}, it is reasonable to expect that
base points from the same $K$-orbit in $X$
determine equivalent invariant algebras.

\begin{definition}
Two invariant algebras $A(X,x_0)$ and $A(X',x_0')$
on $M$ are {\it weakly equivalent}
if $X\cong_G X'$
and there exist $n\in\Aut_G(X)$ and $k\in K$ such that $x_0=n*kx_0'$.
\end{definition}

An invariant algebra $A$ on $M$ may be regarded as an invariant algebra
$\tilde A$ on $K$
such that every element $f\in\tilde A$ is fixed by right $L$-multiplication.
Two such subalgebras $A_1$ and $A_2$ are weakly equivalent if $A_1$ may be
shifted to $A_2$ by the map $R(k): f(x)\to f(xk)$ for some $k\in K$.

Clearly, equivalent invariant algebras are weakly equivalent,
but the converse is not always true.
One may suppose that $x_0=kx_0'$
($\Aut_G(X)$-action does not change the subalgebra).
Consider the subgroups $L_1=K_{x_0}$, $L_2=K_{x_0'}$, and the map
$\phi:K/L\to X$, $\phi(eL)=x_0$. Denote by $\Aut(X,x_0)$ the subgroup
of $\Aut_G(X)$ that preserves $Kx_0$.
(In fact, $\Aut(X,x_0)\subset N_K(L_1)/L_1$.)

\begin{definition}
A closed subgroup $L\subset K$ is {\it an A-subgroup} if any
two weakly equivalent finitely generated invariant algebras on $M=K/L$
are equivalent.
\end{definition}

\begin{proposition} \label{pr1}
A subgroup $L\subset K$ is an A-subgroup
if and only if for any affine embedding $G/F\hookrightarrow X$, $H\subset F$,
and any base point $x_0\in (G/F)^H$ one has
$\Aut(X,x_0)\phi((K/L)^L)=(Kx_0)^L$.
\end{proposition}

{\bf Proof.} Let $x_0'=kx_0$ be an $L$-fixed point.
The equivalence of invariant algebras $\A(X,x_0)$
and $\A(X,x_0')$ means that there is an element $n\in N_K(L)$ such that
$\A(X,nx_0)=\A(X,x_0')$, i.e., $nx_0$ and $x_0'$ are in the same
$\Aut_G(X)$-orbit. If $m\in\Aut_G(X)$ and $m*nx_0=x_0'$, then $m\in\Aut(X,x_0)$.
But the set of points $nx_0$, $n\in N_K(L)$,
coincides with $\phi((K/L)^L)$. $\diamondsuit$

If for any $L\subseteq L_1$ the natural map $(K/L)^L\to (K/L_1)^L$
is surjective, then $L$ is an A-subgroup. In particular, the unit subgroup
and any maximal subgroup in $K$ are A-subgroups.

\begin{corollary}
If $L$ is an A-subgroup, two subgroups $L_1$ and $L_2$ contain $L$
and are $K$-conjugate, then they are $N_K(L)$-conjugate.
\end{corollary}

{\bf Proof.} On $K/L_1$ any point fixed by $L$ has the form
$m*nL_1$, where $m\in N_K(L_1)$ and $n\in N_K(L)$. In particular,
for $L_2=kL_1k^{-1}$, $k\in K$, one has
$kL_1=m*nL_1$ and $L_2=nm^{-1}L_1mn^{-1}=nL_1n^{-1}$. \ $\diamondsuit$

\begin{example}
Put $K=SU(5)$, $L=\{e\}\times\{e\}\times\{e\}\times SU(2)$,
$L_1=SU(2)\times SU(3)$, $L_2=SU(3)\times SU(2)$ as shown on the picture.
Here $L_1$ and $L_2$ are $K$-conjugate, contain $L$, but are not
$N_K(L)$-conjugate. This proves that $L$ is not an A-subgroup.

\bigskip

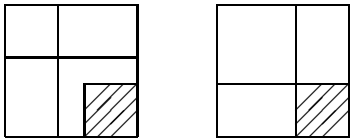
\begin{picture}(200, 50)
\put(10,0){\line(1,0){50}}
\put(60,0){\line(0,1){50}}
\put(10,0){\line(0,1){50}}
\put(10,50){\line(1,0){50}}
\put(90,0){\line(1,0){50}}
\put(140,0){\line(0,1){50}}
\put(90,0){\line(0,1){50}}
\put(10,30){\line(1,0){50}}
\put(30,0){\line(0,1){50}}
\put(40,0){\line(0,1){20}}
\put(40,20){\line(1,0){20}}
\put(90,50){\line(1,0){50}}
\put(40,15){\line(1,1){5}}
\put(40,10){\line(1,1){10}}
\put(40,5){\line(1,1){15}}
\put(40,0){\line(1,1){20}}
\put(45,0){\line(1,1){15}}
\put(50,0){\line(1,1){10}}
\put(55,0){\line(1,1){5}}
\put(120,0){\line(0,1){50}}
\put(90,20){\line(1,0){50}}
\put(120,15){\line(1,1){5}}
\put(120,10){\line(1,1){10}}
\put(120,5){\line(1,1){15}}
\put(120,0){\line(1,1){20}}
\put(125,0){\line(1,1){15}}
\put(130,0){\line(1,1){10}}
\put(135,0){\line(1,1){5}}
\end{picture}
\bigskip
\end{example}


\subsection{Some classes of invariant algebras}\label{vot}

The results of subsection~\ref{ndss} and Theorem~\ref{t1} implies

\begin{proposition}\cite{latdis}\label{propr}
An invariant algebra $A=A(X,x_0)$ is self-conjugate if and only if
$X=Gx_0$ and $G_{x_0}$ is the complexification of $K_{x_0}$.
\end{proposition}

{\bf Remark.} There is one more characterization of this class of $K$-orbits
obtained by V.~M.~Gichev and I.~A.~Latypov. Consider any $G$-equivariant
embedding of $X$ into a $G$-module $V$.
Then the conditions of Proposition~\ref{propr}
are equivalent to the polynomial convexity of the orbit $Kx_0$ in $V$,
see~\cite{gl} for details.

\smallskip

The following theorem due to I.~A.~Latypov may be regarded as
a variant of Luna's theorem (see~\ref{subsecLu}) for compact groups.

\smallskip

\begin{theorem}~\cite{lat} \label{thlat}
Any invariant algebra on $M$ is self-conjugate if and only if the group
$N=N_K(L)/L$ is finite.
\end{theorem}

In this case any invariant algebra on $M$ is finitely generated.
It follows from the results of Section~\ref{secap2} that any invariant
algebra on $M$ is finitely generated if and only if either $N$ is finite
or $K=U(1)$. (Here we assume that the action $K:M$ is effective.)

Now we introduce a class of invariant algebras, which are in some sense
opposite to self-conjugate algebras.

\begin{definition}
An invariant algebra $A$ is said to be {\it
antisymmetric}, if the set $\{ f\in A \mid \overline{f}\in A \}$
coincides with the set of constant functions.
\end{definition}

It is easy to see that antisymmetry is equivalent to any of the following
conditions:

1) any real-valued function in $A$ is a constant;

2) $A$ contains no non-trivial self-conjugate invariant subalgebra.

Hence an invariant algebra $A=A(X,x_0)$ is antisymmetric if and only
if there exists no $G$-equivariant map $\phi: X\to G/H'$, where $G/H'$
is an affine homogeneous space of positive dimension and $G_{\phi(x_0)}$
is the complexification of $K_{\phi(x_0)}$. In particular, if $X$ contains a
$G$-fixed point, then $A(X,x_0)$ is antisymmetric.

\begin{example}
Let $K=SU(2)$, $G=SL(2)$, and $L=H=\{e\}$.
Consider $X=SL(2)/T$. Any point $x_0\in X$ may be regarded as a base point
for some invariant algebra $A(X,x_0)$ on $M=K$. If the stabilizer of $x_0$
contains a torus from $K$, then $A(X,x_0)$ is self-conjugate, and
any two such invariant algebras are equivalent. Other base points determine
antisymmetric algebras: we obtain a 1-parameter family
of mutually non-equivalent
antisymmetric invariant algebras on $SU(2)$.
In particular, this example shows that the property "$A(X,x_0)$
separates points on $M$" depends on the choice of the base point $x_0$ on $X$.
For more information on invariant algebras on $SU(2)$, see~\cite{lat2}.
\end{example}

Finally consider one more natural class of invariant algebras.

\begin{definition}
An invariant algebra $A$ on $M$ is called {\it a Dirichlet algebra} if
the real parts of functions from $A$ are uniformly dense in the algebra of
real-valued continuous functions on $M$.
\end{definition}

Any Dirichlet algebra separates points on $M$, but the converse is not true.
Some results on Dirichlet invariant algebras on compact groups can be found
in~\cite{rid}. In particular, it is proved there that there exists a
biinvariant antisymmetric Dirichlet algebra on $K$ if and only if $K$
is connected and commutative.
It would be interesting to characterize Dirichlet algebras
$A(X,x_0)$ in terms of affine embeddings.


\subsection{Biinvariant algebras and invariant algebras on spheres.}

A {\it biinvariant} algebra on $K$ is a uniformly closed subalgebra
with unit in $C(K)$ invariant with respect to both left and right translations
(here $M=(K\times K)/\Delta(K)$).

Suppose that $F$ is a subgroup in $G\times G$ containing $\Delta(G)$. Then
the subgroup $F_0=\{g\in G \mid (g,e)\in F \}$ is normal in $G$. This
shows that $F$ is the preimage of $\Delta(\tilde G)$
for the homomorphism $G\times G\to \tilde G\times\tilde G$,
where $\tilde G=G/F_0$. Moreover, $\Delta(G)$-fixed points
in $(\tilde G\times\tilde G)/\Delta(\tilde G)$ correspond to central elements
of $\tilde G$. These elements form an orbit of the center $Z(\tilde G)$,
and $Z(\tilde G)$ acts $(\tilde G\times\tilde G)$-equivariantly on any
affine embedding of $(\tilde G\times \tilde G)/\Delta(\tilde G)$.
Hence different
base points on such embeddings define the same invariant algebras. An affine
embedding of the space $(\tilde G\times\tilde G)/\Delta(\tilde G)$ is nothing
else but an algebraic monoid $\tilde S$ with $G(\tilde S)=\tilde G$
(Proposition~\ref{prmon}).

Let us summarize all these observations in the following one-to-one
correspondences (all biinvariant algebras are supposed to be finitely generated):

\smallskip

$\bullet$\ $\{$ self-conjugate biinvariant algebras on $K$ $\}$
$\Longleftrightarrow$ $\{$ quotient groups $\tilde G$ of the group $G$ $\};$

\smallskip

$\bullet$\ $\{$ biinvariant algebras on $K$ $\}$ $\Longleftrightarrow$
$\{$ algebraic monoids $\tilde S$ with $G(\tilde S)=\tilde G~\};$

\smallskip

$\bullet$\ $\{$ biinvariant algebras separating points on $K$ $\}$
$\Longleftrightarrow$ $\{$ algebraic monoids $S$ with $G(S)=G$ $\};$

\smallskip

$\bullet$\ $\{$ antisymmetric biinvariant algebras on $K$ $\}$
$\Longleftrightarrow$ $\{$ algebraic monoids $\tilde S$ with zero
and $G(\tilde S)=\tilde G$ $\}.$

\smallskip

To explain the last equivalence, we note that $\tilde S$ has a zero if and
only if the closed $(\tilde G\times\tilde G)$-orbit in $\tilde S$ is a point.
Embeddings with a $G$-fixed point correspond to antisymmetric invariant
algebras (see~\ref{vot}).
If the closed orbit has positive dimension, it is isomorphic to
$(\tilde G_1\times\tilde G_1)/\Delta(\tilde G_1)$ for a non-trivial quotient
$\tilde G_1$ of the group $\tilde G$,
and the corresponding projection (see~\ref{sslice}) determines a non-trivial
self-conjugate subalgebra in our invariant algebra.

Theorem~\ref{thlat} (or Proposition~\ref{prpol}) shows that any biinvariant
algebra on $K$ is self-conjugate if and only if $K$ is semisimple.
This result was proved by R.~Gangolli~\cite{gan} and J.~Wolf~\cite{wol}.

\medskip

Our final remark concerns invariant algebras on spheres $S^n$.
The classification of transitive actions of compact Lie groups
on spheres was obtained by A.~Borel, D.~Montgomery and H.~Samelson
(see~\cite{on}). All corresponding homogeneous spaces
are spherical with a unique exception: there is a transitive action of
the group $Sp(n)=GL(n,\HH)\cap U(2n)$ on $S^{4n-1}$ with stabilizer
$Sp(n-1)$ and the
complexification of $Sp(n)/Sp(n-1)$ is a homogeneous space of complexity
one. (This is the reason why the clasification of invariant algebras on spheres
was not completed in this case only, see~\cite{latdis}.)

The complexification of $Sp(n)/Sp(n-1)$ satisfies
the conditions of Theorem~\ref{main}.
This implies the following general result: the number of radical invariant
ideals in any invariant algebra on a sphere (with respect to any transitive
action) is finite.


\section{Application Two: $G$-algebras with finitely generated invariant
subalgebras}\label{secap2}

\subsection{The reductive case}
In this section by $\A$ we denote a finitely generated $G$-algebra without
zero divisors. Let us introduce three special types of $G$-algebras.

\medskip

{\it Type C.} \ Here $\A$ is a finitely generated domain
of Krull dimension $\kd \A=1$ (i.e.,
the transcendence degree of the quotient field $Q\A$ equals one)
with any (for example, trivial) $G$-action. Such algebras may be
considered as the algebras of regular functions on
irreducible affine curves.

\medskip

{\it Type HV.}\ Let $\lambda$ be a dominant weight of the group $G$
(with respect to some fixed Borel subgroup)
and $V(\lambda)$ be a simple finite-dimensional $G$-module
with highest weight $\lambda$.
Let $\lambda^*$ be the highest weight
of the dual module $V(\lambda)^*$.
Consider a subsemigroup $P$ in the
additive semigroup of non-negative integers (it is automatically
finitely generated), and put
$$
\A(P,\lambda)=\oplus_{p\in P} V(p\lambda).
$$

There exists a unique structure (up to $G$-isomorphism)
of a $G$-algebra on $\A(P,\lambda)$
such that $V(p\lambda)V(m\lambda)=V((p+m)\lambda)$. In fact,
consider the closure $X(\lambda)=\overline{Gv}$ of the orbit of a
highest weight vector $v$ in $V(\lambda^*)$. The algebra
$\kk[X(\lambda)]$ of regular functions on $X(\lambda)$ as a $G$-module
has the isotypic decomposition
$$
 \kk[X(\lambda)]=\oplus_{k\ge 0}\kk[X(\lambda)]_{k\lambda},
$$
any $\kk[X(\lambda)]_{k\lambda}$ is a simple $G$-module, and
$\kk[X(\lambda)]_{k\lambda}\kk[X(\lambda)]_{m\lambda}=
\kk[X(\lambda)]_{(k+m)\lambda}$, see~\ref{sub1.3}. This allows to realize
$\A(P,\lambda)$ as a subalgebra in $\kk[X(\lambda)]$.
The proof of uniqueness of such multiplication is left to the reader.
Further we shall say that the algebra $\A(P,\lambda)$ is an algebra of
type HV.

\medskip

\begin{example}
Let $G=SL(n)$
and $\omega_1,\dots,\omega_{n-1}$ be its fundamental weights.
The natural linear action
$G:\kk^n$ induces an action on regular functions
$$
 G:\A=\kk[x_1,\dots,x_n], \ \ (gf)(v):=f(g^{-1}v).
$$
The homogeneous polynomials of degree $m$ form an (irreducible)
isotypic component corresponding to the weight $m\omega_{n-1}$.
The algebra $\A$ is of type HV with
$\lambda=\omega_{n-1}$ and $P=\ZZ_+$.
The variety
$X(\omega_{n-1})$ is original space $\kk^n$.
\end{example}

\medskip

{\it Type N.}\ Let $H$ be a closed subgroup of $G$ and
$$
\A(H)=\kk[G]^H=\kk[G/H]=\{f\in \kk[G] \mid f(gh)=f(g) \ {\rm for \ any} \
g\in G,\ h\in H\}.
$$
If $H$ is reductive, then $\A(H)$ is finitely generated.
We say that a $G$-algebra $\A$ is of type N if
there exists a reductive subgroup $H\subset G$ with
$|N_G(H)/H|<\infty$ and $\A$ is $G$-isomorphic to $\A(H)$.

\begin{example}
The algebra $\A(T)=\{f\in \kk[G] \mid f(gt)=f(g) \ \rm {for \ any}\ t\in T\}$
is a $G$-algebra of type N with respect to the left $G$-action.
\end{example}

Now we are ready to formulate the main result.

\smallskip

\begin{theorem}~\cite{ar5} \label{tmain}
Let $\A$ be a finitely generated $G$-algebra without zero divisors.
Then any $G$-invariant subalgebra of $\A$ is finitely
generated if and only if $\A$ is an algebra of one of the types
C, HV or N.
\end{theorem}

We start the proof of Theorem~\ref{tmain} with a method to construct
a non-finitely generated subalgebra.
Let $X$ be an irreducible affine algebraic variety and
$Y$ a proper closed irreducible subvariety. Consider
the subalgebra
$$
  \A(X,Y)=\{f\in \kk[X] \mid f(y_1)=f(y_2) \ {\rm for \ any} \
  y_1, y_2 \in Y\}\subset \A=\kk[X].
$$

\begin{proposition}\label{nfg}
The algebra $\A(X,Y)$ is finitely generated if and only if
$Y$ is a point.
\end{proposition}

{\bf Proof.}\ If $Y$ is a point, then $\A(X,Y)=\kk[X]$. Suppose that
$Y$ has positive dimension and
$\I=\I(Y)=\{ f\in \kk[X] \mid f(y)=0 \ {\rm for \ any}\ y\in Y\}$. Then
$\A/\I$ is infinite-dimensional vector space. By the Nakayama
Lemma, we can find $i\in \I$ such that in the local ring of $Y$
the element $i$ is not in $\I^2$. Then for any
$a\in k[X]\setminus \I$ the element $ia$ is in $\I\setminus \I^2$.
Hence space $\I/\I^2$ has infinite dimension.

On the other hand, suppose that $f_1,\dots,f_n$ are generators
of $\A(X,Y)$. Subtracting constants, one may assume that
all $f_i$ are in $\I$. Then
$\dim \A(X,Y)/\I^2\le n+1$, a contradiction. \ $\diamondsuit$

\begin{proposition}\label{ga}
Let $\A$ be a finitely generated domain. Then any subalgebra in
$\A$ is finitely generated if and only if $\kd \A\le 1$.
\end{proposition}

{\bf Proof.}\ If $\kd \A \ge 2$, then the statement follows from the
previous proposition. The case $\kd \A=0$ is obvious. It remains
to prove that if $\kd \A = 1$, then any subalgebra is finitely
generated. By taking the integral closure,
one may suppose that $\A$ is the
algebra of regular functions on a smooth affine curve $C_1$. Let $C$
be the smooth projective curve such that $C_1\cong C\setminus
\{P_1,\dots,P_k\}$. The elements of $\A$ are rational functions on
$C$ that may have poles only at points $P_i$. Let $\B$ be a subalgebra
in $\A$. By induction on $k$,
we may suppose that the subalgebra $\B'\subset \B$
consisting of functions regular at $P_1$ is finitely generated, say
$\B'=\kk[s_1,\dots,s_m]$. (Functions that are regular at any point
$P_i$ are constants.) Let $v(f)$ be the order of the zero/pole of
$f\in \B$ at $P_1$. The set $V=\{v(f),\ f\in \B\}$ is
an additive subsemigroup
of integers. Such a subsemigroup is finitely generated. Let
$f_1,\dots,f_n$ be elements of $\B$ such that $v(f_i)$ generate $V$.
Then for any $f\in \B$ there exists a polynomial $P(y_1,\dots,y_n)$
with $v(f-P(f_1,\dots,f_n))\ge 0$, thus $f-P(f_1,\dots,f_n)\in
\B'$. This shows that $\B$  is generated by
$f_1,\dots,f_n,s_1,\dots,s_m$. \ $\diamondsuit$

\medskip

Let $\A$ be a finitely generated $G$-algebra with $\kd \A\ge 2$.
Consider the affine variety $X=\Sp \A$.
The action $G:\A$ induces a regular action $G:X$.

Suppose that there exists a proper irreducible closed
invariant subvariety $Y\subset X$ of positive dimension.
Then $\A(X,Y)$ is an invariant subalgebra, which is not finitely
generated. In particular, this is the case if $G$ acts on $X$
without a dense orbit. Hence we may assume that
either

\medskip

(i) the action $G:X$ is transitive, \ or

(ii) $X$ consists of an open orbit and a $G$-fixed point $p$.

\medskip

In case (i), $X=G/H$ and $H$ is reductive. If $G/H$
is not affinely closed then there exists a non-trivial affine embedding
$G/H \hookrightarrow X'$, and the complement in $X$ to the open affine
subset $G/H$ is a union of irreducible divisors.
Let $Y$ be one of these divisors.
The algebra $\A(X',Y)$ is a non-finitely
generated invariant subalgebra in $\kk[X']$
and the inclusion $G/H \hookrightarrow X'$ defines an embedding
$\kk[X']\subset \kk[X]=\A$. On the other hand,

\begin{lemma}
If $X=G/H$ is affinely closed, i.e., $\A$ is
of type N, then any invariant subalgebra in $\A$ is finitely
generated.
\end{lemma}

{\bf Proof.}\ Suppose that there exists an invariant
subalgebra $\B\subset \A$
that is not finitely generated. Let $f_1, f_2,\dots$ be a system
of generators of $\B$. Consider the finitely generated subalgebras
$\B_i=\kk[\langle Gf_1,\dots,Gf_i\rangle]$. Infinitely many of them
are pairwise different.
For the corresponding varieties $X_i:=\Sp \B_i$ one has natural
dominant $G$-morphisms

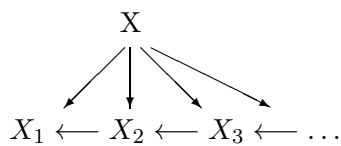
\begin{picture}(200,40)
\put(146,22){X}
\put(148,16){\vector(-1,-1){23}}
\put(154,16){\vector(1,-1){23}}
\put(150,16){\vector(0,-1){23}}
\put(158,16){\vector(2,-1){45}}
\end{picture}
$$
X_1\longleftarrow X_2\longleftarrow X_3\longleftarrow\dots
$$

By Theorem~\ref{afcl}, any $X_i$ is an affine
homogeneous space $G/H_i$, $H\subseteq H_i$. The infinite sequence of
algebraic subgroups
$$
H_1\supset H_2\supset H_3\supset\dots
$$
leads to contradiction.  $\diamondsuit$

\begin{remark}
As is obvious from what has been said any invariant subalgebra in the algebra
$\A(H)$ of type N has the form $\A(H')$, $H\subseteq H'\subseteq G$ and also
has type N. Algebras
of type N can be characterized by the following equivalent properties:

(P1) \  any invariant subalgebra contains no proper invariant ideals;

(P2) \ the algebra contains no proper invariant ideals and the group
of equivariant automorphisms is finite.
\end{remark}

Now consider case (ii).
Let us recall the following theorem due to F.~Bogomolov.

\begin{theorem}\label{tbogo} (~\cite{bo}, see also \cite[Th.7.6]{gr})
Let $X$ be an irreducible affine variety with a non-trivial $G$-action and
with a unique closed orbit, which is a $G$-fixed point.
Then there exists a $G$-equivariant surjective morphism
$\phi: X\to X(\mu)$ for some dominant weight $\mu\ne 0$.
\end{theorem}

In our case the preimage $\phi^{-1}(0)$ is the point $p$, and
thus all fibers of $\phi$ are finite. This shows that $X$ is a
spherical variety of rank one (see~\cite{br} for definitions), i.e.,
$$
\kk[X]=\oplus_{m\ge 0}\kk[X]_{m\lambda},
$$
where $\kk[X]_{m\lambda}$ is either zero or irreducible, and
$\mu=k\lambda$ for some $k>0$. On the other hand, the stabilizer
of any point on $X(\mu)$ contains a maximal unipotent subgroup
of $G$, and the same is true for $X$. By Theorem~\ref{teospec}, this implies
$\kk[X]_{m_1\lambda}\kk[X]_{m_2\lambda}=\kk[X]_{(m_1+m_2)\lambda}$.
Hence $\A=\kk[X]$ is an algebra of type HV.

Conversely, any subalgebra of the $\A(P,\lambda)$ is finitely generated
because it corresponds to some subsemigroup
$P'\subset P$ and $P'$ is finitely generated.\ $\diamondsuit$


\subsection{The non-reductive case}

Let us classify affine $G$-algebras with finitely generated invariant
subalgebras for a non-reductive affine group $G$ with the Levi
decomposition $G=LG^u$. Surprisingly, the result
in this case is simpler than in the reductive case.

In the previous subsection we assumed that a $G$-algebra $\A$ has
no zero divisors. In fact, this restriction is inessential.

\begin{lemma} \cite{ate} \label{lr}
Let ${\rm rad}(\A)$ be the ideal of all nilpotents in $\A$.
The following conditions are equivalent:

{\rm (1)}\ any $G$-invariant subalgebra in $\A$ is finitely generated;

{\rm (2)}\ any $G$-invariant subalgebra in $\A/{\rm rad}(\A)$
is finitely generated and $\dim {\rm rad}(\A)<\infty$.
\end{lemma}

{\bf Proof.} \ Any finite-dimensional subspace in
${\rm rad}(\A)$ generates a finite-dimensional subalgebra in $\A$.
Hence if
$\dim {\rm rad}(\A)=\infty$, then the subalgebra generated by ${\rm rad}(\A)$
is not finitely generated. On the other hand, the preimage in
$\A$ of any non-finitely generated subalgebra in
$\A/{\rm rad}(\A)$ is not finitely generated.

Conversely, assume that (2) holds.
Then any subalgebra in $\A$ is generated by elements whose images
generate the image of this subalgebra in
$\A/{\rm rad}(\A)$, and by a basis of the radical of the subalgebra.
\ $\diamondsuit$

\smallskip

If $\A$ contains non-nilpotent zerodivisors, then the proof of
Theorem~\ref{tmain} goes with small technical modifications, see~\cite{ate}.
The same proof also goes well for a non-reductive $G$.
The only difference is that case HV is excluded by the result of V.~L.~Popov.

\begin{proposition}\cite[Th.3]{po}
If $G$ acts on an affine variety $X$
with an open orbit, and

(1)\ the induced action $G^u:X$ is non-trivial;

(2)\ the complement to the open $G$-orbit in $X$ does not contain a component of
a positive dimension,

then the action $G:X$ is transitive.
\end{proposition}

These arguments prove

\begin{theorem}\cite[Th.3]{ate}
Let $\A$ be a $G$-algebra without nilpotents with the non-trivial
induced $G^u$-action. The following conditions are equivalent:

\smallskip

{\rm (1)}\ any $G$-invariant subalgebra in $\A$ is finitely generated;

{\rm (2)}\ any $G$-invariant subalgebra in $\A$ does not contain
non-trivial $G$-invariant ideals;

{\rm (3)}\ any $L$-invariant subalgebra in $\A^{G^u}$ does not contain
non-trivial $L$-invariant ideals;

{\rm (4)}\ $\A=\kk[G/H]$, where $G/H$ is an affinely closed
homogeneous space.
\end{theorem}


\vspace{0.7cm}

Department of Higher Algebra, Faculty of Mechanics and Mathematics,
Moscow State University, GSP-2, Moscow, 119992, Russia.

\medskip

{\it E-mail address}:\ arjantse@mccme.ru

\medskip

\noindent {\it URL:}
\texttt{http://mech.math.msu.su/department/algebra/staff/arzhan.htm}

\end{document}